\documentclass[10pt]{article}
\usepackage{bbm}
\usepackage{amsmath, amsthm, amssymb}
\usepackage{color}
\usepackage{enumerate}

\newtheorem*{acknowledgements}{Acknowledgements}

\newtheorem*{quec}{Conjecture (QUE for $X={\rm SL}_2 (\mathbbm{Z})\backslash \mathbbm{H}$)}

\newtheorem*{lemma1.3}{Lemma 1.3}
\newtheorem*{theorem1}{Theorem 1}
\newtheorem*{Remark1}{Remark 1}
\newtheorem*{Remark2}{Remark 2}
\newtheorem*{Remark3}{Remark 3}
\newtheorem*{Remark4}{Remark 4}

\newtheorem*{Remark2.1}{Remark}
\newtheorem*{theoremA}{Theorem 2.2}
\newtheorem*{Remark2.3}{Remark}
\newtheorem*{theorem2.3}{Theorem 2.3}

\newtheorem*{lemma3a}{Lemma 3a}
\newtheorem*{lemma3b}{Lemma 3b}
\newtheorem*{lemma3c}{Lemma 3c}
\newtheorem*{lemma3d}{Lemma 3d}

\newtheorem*{lemma4.1}{Lemma 4.1}
\newtheorem*{lemma4.2}{Lemma 4.2}

\newtheorem*{lemma5.1}{Lemma 5.1}
\newtheorem*{lemma5.2}{Lemma 5.2}
\newtheorem*{corollary5.2}{Corollary 5.2}

\newtheorem*{lemmaA.1}{Lemma A.1}
\newtheorem*{lemmaIwaniec}{Lemma A.2a}
\newtheorem*{lemmaA.2b}{Lemma A.2b}
\newtheorem*{theoremA.5}{Theorem A.5}

\newtheorem*{lemmaB.1}{Lemma B.1}
\newtheorem*{corollaryB.2}{Corollary B.2}
\newtheorem*{lemmaB.3}{Lemma B.3}
\newtheorem*{corollaryB.4}{Corollary B.4}
\newtheorem*{theoremB.5}{Theorem B.5}

\begin{document}
\title{A Sieve Method for Shifted Convolution Sums} \author{Roman Holowinsky\thanks{This material is based upon work supported by the National Science Foundation under agreement No. DMS-0111298 and No. DMS-03-01168 and also partially supported by NSERC.  Any opinions, findings and conclusions or recommendations expressed in this material are those of the author and do not necessarily reflect the views of the National Science Foundation.}} \maketitle
\begin{abstract}
We study the average size of shifted convolution summation terms related to the problem of Quantum Unique Ergodicity on ${\rm SL}_2 (\mathbbm{Z})\backslash \mathbbm{H}$.  Establishing an upper-bound sieve method for handling such sums, we achieve an unconditional result which suggests that the average size of the summation terms should be sufficient in application to Quantum Unique Ergodicity.  In other words, cancellations among the summation terms, although welcomed, may not be required. Furthermore, the sieve method may be applied to shifted sums of other multiplicative functions with similar results under suitable conditions.
\end{abstract}

\section{Introduction}
The work presented in this paper focuses on the analysis of the \textit{shifted convolution sums}
\begin{equation}
\sum_{n\leqslant x} \lambda(n) \overline{\lambda(n+\ell)} \label{SCSnoabs}
\end{equation}
where $\ell$ is a fixed non-zero integer and $\lambda$ is a multiplicative function. In particular, we will study such sums when the $\lambda(n)$ are Hecke eigenvalues of Hecke-Maass cusp forms. Obtaining sufficient upper bounds for these types of sums is a key part in the Number Theoretic approach to the problem of Quantum Unique Ergodicity(QUE).  

Instead of analyzing the shifted sums as written in (\ref{SCSnoabs}), we will investigate the behavior of these sums when one ignores possible cancellations among the summation terms by taking absolute values
\begin{equation}
 \sum_{n\leqslant x} | \lambda ( n ) \lambda ( n + \ell ) |. \label{SCSabs}
\end{equation}
Although this may seem wasteful, positivity of the summation terms allows us to proceed with an upper-bound sieve application in mind.  For small fixed $\ell \neq 0$, we benefit from the fact that $n$ and $n+\ell$ have few common factors.  Considering the prime factorizations of $n$ and $n+\ell$ will allow us to in essence remove this additive dependence and demonstrate that obtaining non-trivial upper bounds for the shifted sums (\ref{SCSabs}) will be a direct result of the fact that $|\lambda(n)|<1$ on average.  More specifically, for a given Hecke-Maass cusp form $u$ with Hecke eigenvalues $\lambda$, 
\begin{equation}
\sum_{n\leqslant x} \frac{|\lambda ( n )|}{n} \ll_u (\log x)^{1-\alpha} \label{singleintro}
\end{equation}
for some $\alpha>0$. This type of bound may be established unconditionally (see \S 4.1) by the holomorphicity and positivity at the point $s=1$ of the symmetric square, fourth and sixth power $L$-functions associated with the form $u$.   The upper-bound sieve application will then permit a saving of approximately $(\log x)^\alpha$ for each Hecke eigenvalue factor in order to achieve almost twice as much saving for the shifted sum, i.e. for $\alpha$ the same as in (\ref{singleintro}) and any $\varepsilon>0$ we have
\begin{equation}
\sum_{n\leqslant x} |\lambda(n)\lambda(n+\ell)| \ll_{u,\ell} \frac{x}{(\log x)^{2 \alpha-\varepsilon}}. \label{doubleintro}
\end{equation}
We state this result as our main theorem.
\begin{theorem1}
For a Hecke-Maass cusp form $u$ with Hecke eigenvalues $\lambda ( n )$, a fixed
integer $\ell \neq 0$ and $x\geqslant 2$ we have
  \begin{equation} \sum_{n \leqslant x} | \lambda ( n ) \lambda ( n + \ell ) | \ll_{u,\ell} 
\frac{x}{( \log x )^{\delta}}\label{THM1}   
\end{equation}
  for some absolute positive constant $\delta$. 
\end{theorem1}
\begin{Remark1}
\textnormal{It does not seem as though the method presented in this paper would permit us to obtain a full integral power saving of $\log x$, i.e. $\delta=1$.  We will instead show that (\ref{THM1}) certainly holds with $\delta=1/7$ and this may easily be improved to $\delta \sim 1/5$. However, we do not expect the method in this paper to produce better than $\delta=2(1-8/(3\pi))\sim 0.30235$ for reasons related to the Sato-Tate conjecture.}
\end{Remark1}
\begin{Remark2}\textnormal{One should notice that in each of the bounds (\ref{singleintro}), (\ref{doubleintro}) and (\ref{THM1}), the implied constant depends on the form $u$.  In the problem of QUE, we shall see that this dependence can not be ignored as it will be essential for us to maintain uniform bounds as $u$ varies.}
\end{Remark2}
\begin{Remark3}
\textnormal{The implied constant in (\ref{THM1}) depends on $\ell$ and is bounded by a small power of $\ell$. Therefore, in order to use Theorem 1 successfully in applications, it is necessary to have control over the size of the non-zero shift $\ell$. In Appendix A, we show that one may reduce the problem of QUE to studying shifted sums of type (\ref{SCSabs}) with the size of the non-zero shifts $\ell$ growing arbitrarily slowly relative to the size of the spectral parameter.}
\end{Remark3}
\begin{Remark4}
\textnormal{If one wishes to work under the assumption of the Ramanujan-Petersson(RP) conjecture, arguments may then be simplified while still obtaining Theorem 1 with $\delta \sim 1/9$. For example, the analysis in \S3 involving twists of Rankin-Selberg L-functions may be entirely avoided.  We leave the details of working under the assumption of RP for a separate paper addressing the analogous problem about shifted sums of Hecke eigenvalues for holomorphic cusp forms, as RP is known in this case.  Furthermore, knowledge of RP allows one to consider sums with multiple shifts obtaining similar results.}   
\end{Remark4}

\begin{acknowledgements}
\textnormal{
Many thanks go to everyone involved in the development of this work.  In particular, I thank my thesis advisor Professor Henryk Iwaniec for his influence, expert guidance and devotion to Analytic Number Theory and his students.  I thank Professor Peter Sarnak and Professor Wenzhi Luo for taking the time to meet on several occasions and provide useful suggestions.  Thanks go also to the referees of this paper for the time they took to provide detailed constructive reports. Finally, to the institutions which made it possible to meet with everyone and work on this project, Rutgers University, the Institute for Advanced Study, the University of Toronto and the CRM in Barcelona, where the ideas for this work were first developed, thank you. 
}
\end{acknowledgements}

\subsection{Automorphic forms}
Let $\mathbbm{H}$ be the upper half plane with hyperbolic measure $d \mu z : =
y^{- 2} dxdy$.  Set $\Gamma ={\rm SL}\sb 2(\mathbbm{Z})$ and let $X = \Gamma
\backslash \mathbbm{H}$ be the quotient space with volume
\begin{equation*} \textnormal{Vol} ( X ) : = \int_X d \mu z = \frac{\pi}{3} \end{equation*}
on which we define the Laplace-Beltrami operator
\begin{equation*} \Delta : = y^2 ( \frac{\partial^2}{\partial x^2} +
   \frac{\partial^2}{\partial y^2} ) . \end{equation*}
Let $\mathcal{A}( X )$ be the linear space of \textit{automorphic functions}
\begin{equation*} \mathcal{A}( X ) := \left\{ f :\mathbbm{H} \rightarrow
   \mathbbm{C}\, | \, f ( \gamma z ) = f ( z ) \textnormal{ for all } \gamma \in
   \Gamma \} \right. \end{equation*}
and $\mathcal{A}_s ( X )$ the subspace of \textit{automorphic forms} $f \in
\mathcal{A}( X )$ satisfying
\begin{eqnarray*}
  ( \Delta + s ( 1 - s ) ) f = 0 & \textnormal{with} & s = \frac{1}{2} + it.
\end{eqnarray*}
Denote by $\mathcal{L}( X )$ the Hilbert space of square integrable
automorphic functions with inner product
\begin{equation*} 
< f, g > : = \int_X f ( z ) \bar{g} ( z ) d \mu z. 
\end{equation*}
This space can be decomposed as the closure of two orthogonal subspaces
\begin{equation*} \mathcal{L}( X ) = \tilde{\mathcal{C}} ( X ) \oplus \tilde{\mathcal{E}} ( X
   ) \end{equation*}
where $\mathcal{E}( X )$ is the space of \textit{incomplete Eisenstein} series
\begin{equation*} E ( z| \psi ) := \sum_{\gamma \in \Gamma_{\infty} \backslash \Gamma}
   \psi ( \textnormal{Im } \gamma z ) \end{equation*}
with $\psi$ any smooth compactly supported function on $\mathbbm{R}^+$ and
\begin{equation*}
\Gamma_{\infty} := \left\{ \gamma = \left(\begin{array}{cc}
  1 & \ast\\
  & 1
\end{array}\right) \in \Gamma \right\},
\end{equation*}
and $\mathcal{C}( X )$, its
orthogonal complement, is the space of smooth, bounded automorphic functions
with no zero-th term in their Fourier expansion at $\infty$.  Eigenfunctions of
$\Delta$ in $\mathcal{C}( X )$ are called \textit{cusp forms}.

In addition to the Laplace operator, we have the commuting family of Hecke
operators acting on the space of automorphic functions
\begin{eqnarray*}
  T_n :\mathcal{A}( X ) \longrightarrow \mathcal{A}( X ), & n \geqslant 1. & 
\end{eqnarray*}
The Hecke operators are defined by group operations
\begin{equation*} ( T_n f ) ( z ) = \frac{1}{\sqrt{n}}  \sum_{\tau \in \Gamma \backslash G_n} f (
   \tau z ), \end{equation*}
where 
\begin{equation*}
G_n = \left\{ \left(\begin{array}{cc}
  a & b\\
  c & d
\end{array}\right) : a, b, c, d \in \mathbbm{Z}, ad - bc = n \right\},
\end{equation*} 
and thus commute with $\Delta$.  Eigenvalues $\lambda_f ( n )$ of the Hecke
operators are multiplicative and satisfy the rule
\begin{equation*} \lambda_f ( m ) \lambda_f ( n ) = \sum_{d| ( m, n )} \lambda_f \left(
   \frac{mn}{d^2} \right) . \end{equation*}
Since these Hecke operators commute with one another and are self-adjoint, we
can choose an orthogonal basis $\{ u_j ( z ) \}$ for the space of cusp forms
$\mathcal{C}( X )$, such that for all $n$,
\begin{equation*} T_n u_j ( z ) = \lambda_j ( n ) u_j ( z ) . \end{equation*}
Functions $u$ which are eigenfunctions of $\Delta$ and simultaneous
eigenfunctions for all $T_n$ are called \textit{Hecke-Maass cusp forms}.  For
convenience, we normalize these forms so that
\begin{equation*} < u_j, u_j > = \int_X |u_j ( z ) |^2 d \mu z = 1 . \end{equation*}
With this normalization, a Hecke-Maass cusp form with $\Delta$ eigenvalue $\lambda =
1 / 4 + t_j^2$ will have a Fourier expansion of the type
\begin{eqnarray}
  \nonumber u_j ( z ) & = & \sum_{n \neq 0} a_n ( y ) e ( nx )\\
  & = & \sqrt{y}  \sum_{n \neq 0} \rho_j ( n ) K_{i t_j} ( 2 \pi |n|y ) e ( nx ) \label{0.5}
\end{eqnarray}
where $K$ is the $K$-Bessel function defined in Appendix B and $\rho_j ( n )$ is
proportional to the $n$-th Hecke eigenvalue $\lambda_j ( n )$ of $u_j$
\begin{eqnarray*} 
\rho_j ( n ) = \lambda_j ( n ) \rho_j ( 1 ) & \textnormal{for} & n\geqslant 1. 
\end{eqnarray*}
We then have that either $u_j(z)$ is \textit{even}, when $\rho_j(-n)=\rho_j(n)$, or that $u_j(z)$ is \textit{odd}, when $\rho_j(-n)=-\rho_j(n)$.
\subsection{Quantum Unique Ergodicity}
Let $\{u_j\}$ be an orthonormal basis of Hecke-Maass cusp forms.  Our work is motivated by the following conjecture.
\begin{quec}
For any $f$ in $\mathcal{L}( X )$ we have
\begin{eqnarray}
  < fu_j, u_j > = \frac{3}{\pi} < f, 1 > + o ( 1 ) & \textnormal{as} & j \longrightarrow
  \infty  \label{1}
\end{eqnarray}
\end{quec}
\noindent Note that when $f$ is a constant function then (\ref{1}) is immediate since 
\begin{equation*}
< cu_j, u_j > = c < u_j, u_j > = c
\end{equation*}
for any constant $c$.

This conjecture and several analogues, originally stated in [R-S] for general compact manifolds of negative curvature, have been studied by a variety of authors. The most notable result towards establishing such an equidistribution statement is due to Lindenstrauss [Li], in which he succeeded to establish QUE completely for $X$ a compact arithmetic quotient in $\mathbbm{H}$.  When $X$ is not compact, the conjecture holds by Lindenstrauss's method for cuspidal $f$, but still requires subconvexity estimates for the symmetric square $L$-function of $u_j$ for $f$ an incomplete Eisenstein series. 

If one is looking for effective rates of convergence or studying the analogous equidistribution problem in the holomorphic setting, however, one must currently apply techniques from Analytic Number Theory.  The work of Luo and Sarnak [L-S], for example, when the Hecke-Maass cusp forms $\{u_j\}$ are replaced with holomorphic forms with weight going to infinity, involves reducing the equidistribution problem to the study of shifted convolution sums and the analysis of Kloosterman sums.

The approach we take in this paper will be similar to that of Luo and Sarnak in that our main concern is the analysis of shifted sums of Hecke eigenvalues.  The difference in our work comes from the belief that the average size of the shifted summation terms should be sufficient for QUE.  Therefore, we reduce the problem to the study of the sums
\begin{equation}
\sum_{n\asymp t_j} | \lambda_j(n) \lambda_j(n+\ell)| \label{SCSabswithj}
\end{equation}
for $\ell \neq 0$, instead of the sums 
\begin{equation}
\sum_{n\asymp t_j} \lambda_j(n) \lambda_j(n+\ell), \label{SCSnoabswithj}
\end{equation}
and will apply upper-bound sieve techniques instead of Kloosterman sum analysis. The method of reduction will in some ways be simpler, as expected by taking absolute values, and will also demonstrate that the number of shifts $\ell \neq 0$ which need to be considered will grow arbitrarily slowly relative to the spectral parameter $t_j$. In fact, in Appendix A we show that modulo a ``log-free" convexity bound for the Rankin-Selberg zeta function associated with $u_j$ the problem of QUE reduces to proving 
\begin{equation}
  \sum_{n\asymp t_j} \left|{\lambda}_j ( n ) \lambda_j ( n + \ell )\right| = o ( t_j L(\textnormal{sym}^2 u_j, 1) ) \label{ultimategoal}
\end{equation}
as $t_j\longrightarrow \infty$ for $0 < |\ell|  \leqslant L$ with $L\longrightarrow \infty$ arbitrarily slowly.  

\subsection{QUE obstacles due to the spectral parameter $t_j$}
Recall that our bound in Theorem 1 has an implied constant depending on the Hecke-Maass cusp form in question.  At present, the dependence on the spectral parameter does not permit us to establish (\ref{ultimategoal}).  Starting with the partial sum
\begin{equation*}
\sum_{n \leqslant x} | \lambda_j ( n ) \lambda_j ( n + \ell ) |,
\end{equation*}
the sieve method in this paper ultimately and unconditionally brings one to the bound below from which Theorem 1 will follow (see \S2).
\begin{lemma1.3}
For a Hecke-Maass cusp form $u_j$ with Hecke eigenvalues $\lambda_j( n )$, a fixed
integer $\ell \neq 0$ and $x\geqslant 2$ we have
 \begin{equation} \sum_{n \leqslant x} | \lambda_j ( n ) \lambda_j ( n + \ell ) | \ll_{\ell} x L(\textnormal{sym}^2 u_j, 1) M_j(x) +R_j(x) \label{THM3}   
\end{equation} 
where 
\begin{equation}
M_j(x):=\prod_{p \leqslant z}\left(1-\frac{(1-|\lambda_j(p)|)^2}{p}\right) \label{MTHM3}
\end{equation}
with $z=x^{1/s}$ and $s=c \log \log x$ for some large constant $c$ and $R_j(x)$ is an ``error" term depending on $u_j$ of size
\begin{equation}
R_j(x)=O_j\left(x (\log x)^{-C}\right)\label{RTHM3}
\end{equation}
for any $C>0$.
\end{lemma1.3}

When $x$ is of size $t_j$ as in (\ref{ultimategoal}), showing that $R_j(x)$ is smaller than the main term with respect to $t_j$ currently requires one to assume the Ramanujan-Petersson conjecture as well as subconvexity bounds for Dirichlet character twists of the Rankin-Selberg zeta function associated with the Hecke-Maass cusp form $u_j$.

Assuming uniform bounds on the error term $R_j(x)$ in (\ref{THM3}), it is interesting to note that failing to establish (\ref{ultimategoal}) would therefore be due to a peculiar bias in the Hecke eigenvalues. If for prime $p \leqslant z$ we would have $\lambda_j(p)$ taking values sufficiently close to $\pm 1$, then the main term factor $M_j(x)$ would provide no additional saving, leaving us with 
\begin{equation*} 
\sum_{n \leqslant t_j} | \lambda_j ( n ) \lambda_j ( n + \ell ) | \ll_{\ell} t_j L(\textnormal{sym}^2 u_j, 1)  
\end{equation*}
instead of the desired ``small-o" estimate in (\ref{ultimategoal}). Although such a bias is not expected, present estimates are not strong enough to rule out the possibility.  

In \S4, however, we see that the main term factor $M_j(x)$ may be bounded in terms of information coming from the associated symmetric power $L$-functions. In particular, Lemma 4.1 states that
\begin{equation*}
M_j(x) \ll \left(\frac{L_6(u_j,z)}{L_2(u_j,z) L^2_4(u_j,z) (\log z)^3}\right)^{1/18}
\end{equation*}
with $L_m(u,z)$ the partial Euler product of the m-th symmetric power $L$-function associated with a form $u$ evaluated at the point $s=1$ (see \S 4)
\begin{equation*}
L_{m}(u,z) =  \prod_{p\leqslant z} \prod_{j=0}^m\left(1-\alpha_p^{m-j}\beta_p^{j} p^{-1}\right)^{-1}.
\end{equation*}
Therefore, 
\begin{equation}
\sum_{n\leqslant t_j} \left|{\lambda}_j ( n ) \lambda_j ( n + \ell )\right| \ll_\ell \frac{t_j L(\textnormal{sym}^2 u_j, 1)} {\left(L_2(u_j,z) L^2_4(u_j,z)L^{-1}_6(u_j,z) (\log z)^3\right)^{1/18}}. \label{symfactors}
\end{equation}
Note that the Grand Reimann Hypothesis would have the full symmetric powers satisfying
\begin{equation*}
(\log \log t_j)^{-B_m} \ll L(\textnormal{sym}^m u_j, 1) \ll (\log \log t_j)^{B_m}
\end{equation*}
for some constant $B_m$ depending on the symmetric power $m$, hence we expect to achieve some logarithmic saving uniform in $t_j$ in (\ref{symfactors}) and therefore establish (\ref{ultimategoal}), although an unconditional proof is still out of our reach. It would therefore be interesting to study how partial products of symmetric power $L$-functions behave with respect to one another in the hope of showing
\begin{eqnarray*}
L_2(u_j,z) L^2_4(u_j,z) L^{-1}_6(u_j,z) (\log z)^3 \longrightarrow \infty & \textnormal{as} & t_j\longrightarrow \infty
\end{eqnarray*}
and other similar relations.
 
\newpage
\section{Proof of Theorem 1, the sieve method}
As mentioned in the introduction, we seek to benefit from the fact that $n$ and $n+\ell$ have few common factors when $\ell$ is small.  We first rearrange our shifted sums, which we shall from now on denote as
\begin{equation}
S_{\ell} ( x ) := \sum_{n \leqslant x} | \lambda ( n ) \lambda ( n + \ell ) |,
\label{SCS}
\end{equation}
into an object to which a double upper-bound sieve may be applied.  Note that we have dropped the subscript-$j$ notation for the sake of presentation.  Afterwards, we will demonstrate how our analysis brings one to the bound (\ref{THM3}) in \S 1.3.  Applying results from \S 4 will then prove Theorem 1.
\subsection{Factorization and partitioning}
Let $P(z)$ be defined as
\begin{equation}
P(z):=\prod_{p\leqslant z} p.
\label{P}
\end{equation}
Factoring $n$ and $n + \ell$ uniquely as
\begin{eqnarray}
  n = ab & \textnormal{and} & n + \ell = a_\ell b_\ell, \label{9.2.1}
\end{eqnarray}
such that for every prime $p$ dividing $n(n+\ell)$,
\begin{eqnarray}
 p| a a_\ell  \Rightarrow p \leqslant z & \textnormal{and} & (b b_\ell, P(z))=1 \label{9.2.2}
\end{eqnarray}
with $z=x^{1/s}$ for some $s$ to be chosen later, we partition the sum $S_{\ell} ( x )$ into parts depending on the size of $a$ and $a_\ell$.  The reasoning behind this partition is that we expect to easily treat the parts where $a$ or $a_\ell$ are large, because the number of such $a$ and $a_\ell$ with small prime factors should be relatively small. We denote by $\mathcal{S}_A ( x )$ and $\mathcal{S}_{A_\ell}( x )$ the parts where $a$ and $a_\ell$ are greater than $x^{1/16}$ 
\begin{eqnarray}
  \mathcal{S}_A ( x ) & := & \sum_{\substack{
    n = ab \leqslant x\\
    n + \ell = a_\ell b_\ell\\
    p| aa_\ell  \Rightarrow p \leqslant z\\
    (b b_\ell, P(z))=1\\
    a > x^{1 / 16}
  }} | \lambda ( n ) \lambda ( n + \ell ) | \label{SA}\\
  \mathcal{S}_{A_\ell} ( x ) & := & \sum_{\substack{
    n = ab \leqslant x\\
    n + \ell = a_\ell b_\ell\\
    p| aa_\ell \Rightarrow p \leqslant z\\
    (b b_\ell, P(z))=1\\
    a_\ell > x^{1 / 16}
  }} | \lambda ( n ) \lambda ( n + \ell ) | \label{SAl}
\end{eqnarray}
and in \S 5.1, Lemma 5.1, will show that
\begin{equation}
\mathcal{S}_A(x)+\mathcal{S}_{A_\ell}(x) \ll_{u,\ell} x (\log x)^{-C}
\end{equation}
for any constant $C > 0$.

The sum over the part with both $a$ and $a_\ell$ less than or equal to $x^{1/16}$ we denote by $S^{\star}_{\ell} ( x )$ 
\begin{equation}
S^{\star}_{\ell} ( x ) := \sum_{\substack{
    n = ab \leqslant x\\
    n + \ell = a_\ell b_\ell\\
    p| aa_\ell \Rightarrow p \leqslant z\\
    (b b_\ell, P(z))=1\\
    a, a_\ell \leqslant x^{1 / 16}
  }} | \lambda ( n ) \lambda ( n+\ell )
  | \label{9.3.3}
\end{equation}
so that
\begin{equation}
S_{\ell}(x)\leqslant \mathcal{S}_A(x)+\mathcal{S}_{A_\ell}(x)+S^{\star}_{\ell} (x) \ll_{u,\ell} x (\log x)^{-C} + S^{\star}_{\ell}(x).
\label{4parts}
\end{equation}
Writing $n$ and $n+\ell$ in terms of their factorizations, we are left with evaluating
\begin{equation}
S^{\star}_{\ell}(x)=\sum_{\substack{
     a, a_\ell \leqslant x^{1 / 16}\\
     p|aa_\ell \Rightarrow p\leqslant z
   }} | \lambda ( a ) \lambda ( a_\ell ) |
   \sum_{\substack{
     b \leqslant x / a\\
     a_\ell b_\ell = ab + \ell\\
     ( bb_\ell, P ( z ) ) = 1
   }} |\lambda ( b ) \lambda(b_\ell)|. \label{writeasfactors}
\end{equation}
To help deal with certain co-primality conditions which come up during analysis, we denote the greatest common divisor of $a$ and $a_\ell$ as $v$ and have 
\begin{equation}
S^{\star}_{\ell}(x) \leqslant \sum_{v|\ell} \sum_{\substack{
     a, a_\ell \leqslant x^{1 / 16}\\
     p|aa_{\ell} \Rightarrow p\leqslant z\\
     (a,a_\ell)=1
   }} | \lambda ( v a ) \lambda ( v a_\ell) | \sum_{\substack{
     b \leqslant x / a\\
     a_\ell b_\ell = ab + \ell' \\
     ( bb_\ell, P ( z ) ) = 1
   }} |\lambda ( b ) \lambda(b_\ell)|\label{dividebyv}
\end{equation}
with $\ell'=\ell/v$.  In \S2.3 we will treat the inner most sum in (\ref{dividebyv}) by an application of a double-upper bound sieve which we shall now describe.
\subsection{Double upper-bound sieve heuristics}
Let $\mathcal{A}= ( a_n )$ be a sequence of non-negative numbers on which
we would like to obtain estimates for the sifted sum
\begin{equation*} S ( x, z ) : = \sum_{\substack{
     n \leqslant x\\
     ( n, P ( z ) ) = 1
   }} a_n . \end{equation*}
Here $P ( z )$ is the product of primes less than $z$ which belong to some set
$\mathcal{P}$
\begin{equation*} P ( z ) : = \prod_{\substack{
     p \leqslant z\\
     p \in \mathcal{P}
   }} p. \end{equation*}
Application of an upper-bound sieve with linear sieve weights $\xi_{d}$ of level $D$ 
\begin{equation*}
S(x,z)\leqslant \sum_{\substack{
     d|P(z)\\
     d\leqslant D
   }}\xi_{d} \sum_{\substack{
     n \leqslant x\\
     n \equiv 0 ( \textnormal{mod } d )
   }} a_n 
\end{equation*}
removes the co-primality condition and
allows us to study the partial sums
\begin{equation*} A_d ( x ) : = \sum_{\substack{
     n \leqslant x\\
     n \equiv 0 ( \textnormal{mod } d )
   }} a_n . \end{equation*}
If asymptotics can be established for $A_d ( x )$ when $d|P ( z )$ of the form
\begin{equation*} A_d ( x ) = g ( d ) X + r_d ( x ) \end{equation*}
where $g ( d ) X$ is the expected main term with ``density" $g(d)$ a multiplicative function and $r_d ( x )$ is a sufficiently
small error term, then one would expect to establish a bound
\begin{equation*} 
S ( x, z ) \ll XV 
\end{equation*}
with $X$ approximately equal to
\begin{equation*} A ( x ) := \sum_{n \leqslant x} a_n \end{equation*}
and
\begin{equation*} V = \prod_{p|P ( z )} ( 1 - g ( p ) ) . \end{equation*}
Establishing asymptotics for $A_d(x)$ is necessary as our sieve weights can take positive and negative values.
 
Taking the example of $g ( p ) : = p^{- 1}$ and $P ( z )$ the product of all
primes up to $z$, Merten's formula would then give a saving of $\log z$
\begin{equation*} S ( x, z ) \ll X ( \log z )^{- 1} . \end{equation*}
In general, if a multiplicative function $g$ satisfies $0 \leqslant g ( p ) <
1$ and ``behaves'' like $p^{- 1}$ at primes $p$, then one should see similar
results.

This idea can be extended to a double sum
\begin{equation*} S ( x, z ) := \sum_{\substack{
     n_1 \leqslant x\\
     ( n_1, P ( z ) ) = 1
   }} \sum_{\substack{
     n_2 \leqslant x\\
     ( n_2, P ( z ) ) = 1
   }} a_{n_1, n_2} \end{equation*}
for some sequence of non-negative numbers $\mathcal{A}= ( a_{n_1, n_2} )$. 
In this case, we study the partial sums
\begin{equation*} A_{d_1, d_2} ( x ) := \sum_{\substack{
     n_1 \leqslant x\\
     n_1 \equiv 0 ( \textnormal{mod } d_1 )
   }} \sum_{\substack{
     n_2 \leqslant x\\
     n_2 \equiv 0 ( \textnormal{mod } d_2 )
   }} a_{n_1, n_2} \end{equation*}
and look for asymptotics of the form
\begin{equation} A_{d_1, d_2} ( x ) = g' ( d_1 ) g'' ( d_2 ) X + r_{d_1, d_2} ( x ) \label{DAF} \end{equation}
as before.  Provided that the appropriate conditions are satisfied, we carry
through in the same manner to obtain the bound
\begin{equation} S ( x, z ) \ll XV' V'' \label{DUB} \end{equation}
where
\begin{eqnarray*}
  V' = \prod_{p|P ( z )} ( 1 - g' ( p ) ) & \textnormal{and} & V'' = \prod_{p|P ( z
  )} ( 1 - g'' ( p ) ) .
\end{eqnarray*}
However, if the density functions $g'$ and $g''$ are not independent in terms
of $d_1$ and $d_2$, then establishing (\ref{DUB}) would require some more care. For example, if asymptotics in (\ref{DAF}) were to hold only when $( d_1, d_2 ) = 1$, then upon application of a double upper-bound sieve with sieve weights $\xi'_{d_1}$ and $\xi''_{d_2}$, of level $D'$ and $D''$ respectively, one would arrive at
\begin{equation*} S ( x, z ) \ll X ( G' \ast G'' ) \end{equation*}
where
\begin{equation} G' \ast G'' := \sum_{\substack{d_1 |P ( z )\\ d_1<D'}} \xi'_{d_1} g' ( d_1 )
   \sum_{\substack{
     d_2 |P ( z )\\
     d_2 <D''\\
     ( d_1, d_2 ) = 1
   }} \xi''_{d_2} g'' ( d_2 ) . \label{G'*G''}
\end{equation}
In the appendix of [F-I], Friedlander and Iwaniec treat such a condition. We state their theorem here with our definition of $G'\ast G''$.
\begin{theoremA}
  Let $g'$ and $g''$ be density functions satisfying the linear sieve
  conditions. Let $( \xi' )$ and $( \xi'' )$ be the optimal linear sieve
  weights of level $D'$ and $D''$ respectively. Then we have
  \begin{equation*} G' \ast G'' \ll CV' V''  \end{equation*}
  where
  \begin{eqnarray*}
    V' = \prod_{\substack{
      p|P ( z )\\
      p < D'
    }} ( 1 - g' ( p ) ), & &
    V'' = \prod_{\substack{
      p|P ( z )\\
      p < D''
    }} ( 1 - g'' ( p ) )  
  \end{eqnarray*}
  and
  \begin{equation*} C = \prod_{p|P ( z )} ( 1 + h' h'' ( p ) ) \end{equation*}
  with
  \begin{eqnarray*}
    h' ( p ) = g' ( p ) ( 1 - g' ( p ) )^{- 1}, & h'' ( p ) = g'' ( p ) ( 1 - g''( p ) )^{- 1} . & 
  \end{eqnarray*}
\end{theoremA}
Therefore, if $C$ is uniformly bounded in $z$, the dependence of
$g''$ on $g'$ by the condition $g'' ( d_2 ) = 0$ when $( d_1, d_2 ) \neq 1$ can
be removed and we once again obtain the upper bound (\ref{DUB}). Furthermore, if both density functions behave like $g(p)=p^{-1}$, then we expect to save an extra power of $\log z$ and obtain
\begin{equation*} S ( x, z ) \ll X ( \log z )^{- 2} . \end{equation*}
\subsection{Applying sieve to the inner sum in (\ref{dividebyv})}
In order for our sieve method to work properly in application to the sums seen in (\ref{dividebyv})
\begin{equation*}
\sum_{\substack{
     b \leqslant x / a\\
     a_\ell b_\ell = ab + \ell\\
     ( bb_\ell, P ( z ) ) = 1
   }} |\lambda ( b ) \lambda(b_\ell)|
\end{equation*}
with fixed $(a, a_\ell)=1$, we must choose a sifting sequence $\mathcal{A}$ for which we will be able to establish asymptotics of the form (\ref{DAF}) when considering the relevant partial sums in arithmetic progressions.  To this end, we make use of the trivial bound
\begin{equation*}
|\lambda(b)\lambda(b_\ell)|\leqslant \lambda^2(b)+\lambda^2(b_\ell)
\end{equation*}
and set $\mathcal{A}=(\lambda^2(b))$ for $b,b_\ell$ satisfying $a_\ell b_\ell=ab+\ell$. By symmetry, the sequence $(\lambda^2(b_\ell))$ is treated the same. We may also assume that our sequence is supported on square-free integers $b$, redefining $\mathcal{A}=(\lambda^2(b) \mu^2(b))$, as Lemma 5.2 will demonstrate that the contribution from those $b$ containing square factors will be quite small.  This is done for the purely technical reason that $\lambda^2$ is not a completely multiplicative function.  Furthermore, we loosen the co-primality conditions on $b$ and $b_\ell$ by positivity so that we now sum over $( bb_\ell, P_{6\ell} ( z ) ) = 1$ with
\begin{equation}
P_{6\ell} ( z ) :=\prod_{\substack{p\leqslant z\\ p \nmid 6\ell}} p.\label{Ph}
\end{equation}
Finally, we take advantage of positivity once more by defining a multiplicative function $\eta$ such that
\begin{equation}
\eta(p^\alpha):=\left\{\begin{array}{ll}
\lambda^2(p) & \textnormal{if }  p \mid 6 \\
\lambda^{2\alpha}(p) & \textnormal{if }  p \nmid 6
\end{array}\right. \label{eta}
\end{equation}
for prime powers $p^\alpha$. On square-free $b$, we see that $\eta(b)=\lambda^2(b) \mu^2(b) = \lambda^2(b)$. For all other $b$, we have $\eta(b) \geqslant \lambda^2(b) \mu^2(b) = 0$  and can redefine our sifting sequence for the final time as $\mathcal{A}=(\eta(b))$ for $b, b_\ell$ satisfying $a_\ell b_\ell=ab+\ell$. We are therefore left with applying the double upper-bound sieve to the sum
\begin{equation}
S_{a,a_\ell}(x,z):= \sum_{\substack{
     b \leqslant x / a\\
     a_\ell b_\ell = ab + \ell\\
     ( bb_\ell, P_{6\ell} ( z ) ) = 1
   }} \eta(b) \label{sievingsum}.
\end{equation}

For $(a d, a_\ell d_\ell)=1$, we seek asymptotics for the partial sums in arithmetic progression
\begin{equation} A_{d, d_\ell} ( x ) := \sum_{\substack{
     c \leqslant x/ad\\
     a_\ell d_\ell c_\ell= a d c + \ell
   }} \eta ( d c )
\label{3sum}
\end{equation}
to be of the form
\begin{equation} A_{d, d_\ell} ( x ) = g' ( d ) g'' ( d_\ell ) X + r_{d, d_\ell} ( x ) 
\label{3goal}
\end{equation}
where $g' ( d ) g'' ( d_\ell ) X$ is the expected main term with $g'$ and $g''$ multiplicative functions and the error term
$r_{d, d_\ell} ( x )$ is sufficiently small.  Here $X$ is approximately equal to 
\begin{equation*}
A(x)=\sum_{\substack{ab\leqslant x\\ a_\ell b_\ell = a b + \ell}} \eta(b).
\end{equation*}
In \S 3 we prove the following asymptotic using standard analysis of twisted Rankin-Selberg $L$-functions.
\begin{theorem2.3}
Fix an integer $\ell \neq 0$. Fix integers $a, a_\ell, d, d_\ell \leqslant x^{1/16}$
such that $(a d, a_\ell d_\ell) = 1$ and $d, d_\ell |P_{6\ell} ( z )$.
For $A_{d,d_\ell}(x)$ as defined in (\ref{3sum}) we have
\begin{equation} 
A_{d, d_\ell} ( x ) = \frac{\theta}{\zeta(2)} L(\textnormal{sym}^2 u,1) \frac{\lambda^2(d)}{\varphi(a_\ell d_\ell)} \frac{x}{a d} + O_u(x^{13/16+\varepsilon}) \label{THM2.3}
\end{equation}
for any $\varepsilon>0$ where
\begin{equation*}
\theta=\gamma_u \prod_{\substack{p\mid a_\ell \\ p\mid 6}} \left(1+\frac{\lambda^2(p)}{p}\right)^{-1}\prod_{\substack{p\mid a_\ell d_\ell \\ p\nmid 6}} \left(1-\frac{\lambda^2(p)}{p}\right)
\end{equation*}
with constant $\gamma_u$ satisfying
\begin{equation*}
\frac{3}{5\pi^2} < \gamma_u < 15. 
\end{equation*}
\end{theorem2.3}

To have our asymptotic (\ref{THM2.3}) be of the form (\ref{3goal}), we set
\begin{equation}
h(q):=\prod_{p\mid q} (1-\lambda^2(p) p^{-1})
\end{equation}
and define our density functions $g'$ and $g''$ to be
\begin{equation}
g'(p):=\left\{\begin{array}{ll}
    \lambda^2(p) p^{-1}  & \textnormal{if } p \nmid a_\ell\\
    0 & \textnormal{if } p \mid a_\ell
  \end{array}\right. \label{multfunc3g'}
\end{equation}
and
\begin{equation} 
g''(p):=\left\{\begin{array}{ll}
    h(p) \varphi(p)^{-1} & \textnormal{if } p \nmid a a_\ell\\
    p^{-1} & \textnormal{if } p \nmid a, p\mid a_\ell\\
    0 & \textnormal{if } p\mid a 
  \end{array}\right..
\label{multfunc3g''}
\end{equation}
Observe that 
\begin{equation}
g''(p) = \frac{1}{p}+O(p^{-3/2}) \label{g''close}
\end{equation}
for $p\nmid a$. With these definitions, we bring our asymptotic to the desired form (\ref{3goal}) with main term 
\begin{equation}
X=\gamma_u L(\textnormal{sym}^2 u,1) \frac{x}{a} \frac{h(a_\ell)}{\varphi(a_\ell)}
\label{X3}
\end{equation}
and error term
\begin{equation}
r_{d,d_\ell} (x) = O_u(x^{13/16+\varepsilon}). \label{mainerrors}
\end{equation}

Choosing linear sieve weights $\xi'_{d}$ and $\xi''_{d_\ell}$ of level $D=x^{1/64}$, the contribution from the error term $r_{d,d_\ell}(x)$ will be negligable with respect to $x$ as 
\begin{equation*}
\sum_{\substack{d|P_{6\ell} (z) \\ d \leqslant D}}|\xi_{d}'|\sum_{\substack{d_{\ell}|P_{6\ell} (z)\\ d_\ell \leqslant D\\}}|\xi_{d_\ell}''| \ll (D \log^2 D)^2.
\end{equation*}
This is because linear sieve properties give $|\xi '|, |\xi ''|\leqslant \tau_3$, where $\tau_3(n)$ is the number of ways to write $n$ as the product of three integers.  Therefore, we will have the bound
\begin{equation}
S_{a,a_\ell}(x,z) \leqslant X \sum_{\substack{
     d |P_{6\ell} ( z )\\
     d \leqslant D\\
     (d, a_\ell)=1
      }} \xi_{d}' g'(d)
   \sum_{\substack{
     d_\ell |P_{6\ell} ( z )\\
     d_\ell \leqslant D\\
     (d_\ell, a)=1\\
     ( d, d_\ell ) = 1
   }} \xi''_{d_\ell} g''(d_\ell)  + O_u\left(\frac{x^{63/64}}{a a_\ell}\right)\label{mainupper}
\end{equation}
so that an application of Theorem 2.2 gives 
\begin{equation}
S_{a,a_\ell}(x,z) \ll X C V' V''\label{oursievebound}
\end{equation}
with
\begin{eqnarray*}
C V' V'' & = & \prod_{\substack{p\leqslant z \\ p \nmid 6\ell}}\left(1+\frac{g'(p) g''(p)}{(1-g'(p))(1-g''(p))}\right)\prod_{\substack{p\leqslant z \\ p \nmid 6\ell}}(1-g'(p))(1-g''(p))\\
& \ll_\ell & \prod_{p\leqslant z}\left(1+\frac{\lambda^2(p)}{p^2}\right)\prod_{p\leqslant z}(1-g'(p))(1-g''(p))\\
& \ll_\ell & \prod_{p\leqslant z}(1-g'(p))(1-g''(p)).
\end{eqnarray*} 
\subsection{Establishing Lemma 1.3 and concluding Theorem 1}
With our main term $X$ as in (\ref{X3}) we have,
\begin{equation*}
\frac{h(a_\ell)}{\varphi(a_\ell)} \prod_{p\leqslant z}(1-g'(p)) \ll \frac{1}{\varphi(a_\ell)}\prod_{p\leqslant z}\left(1-\frac{\lambda^2(p)}{p}\right).
\end{equation*}
By (\ref{g''close}) we have,
\begin{equation*}
\frac{1}{a}\prod_{p\leqslant z}(1-g''(p)) \ll \frac{1}{\varphi(a)}\prod_{p\leqslant z}\left(1-\frac{1}{p}\right).
\end{equation*}
Therefore,
\begin{eqnarray*}
S^{\star}_{\ell}(x) & \ll_\ell & x L(\textnormal{sym}^2 u,1) \sum_{\substack{
     a, a_\ell \leqslant x^{1 / 16}\\
     p|aa_\ell \Rightarrow p\leqslant z
   }} \frac{| \lambda ( a ) \lambda ( a_\ell ) |}{ \varphi(a)  \varphi(a_\ell)} \prod_{p\leqslant z}\left(1-\frac{\lambda^2(p)}{p}\right)\left(1-\frac{1}{p}\right)\\
& \ll_\ell & x L(\textnormal{sym}^2 u,1) \prod_{p\leqslant z} \left(1+\frac{|\lambda(p)|}{p}\right)^2 \left(1-\frac{\lambda^2(p)}{p}\right)\left(1-\frac{1}{p}\right)\\
 & \ll_\ell & x L(\textnormal{sym}^2 u,1) M_u(x)
\end{eqnarray*} 
with $M_u(x)$ as in (\ref{MTHM3}). Combining this with (\ref{4parts}) concludes the proof of Lemma 1.3.
\begin{Remark2.3}
\textnormal{We see that the error term $R_u(x)$ in Lemma 1.3 comes not only from the sums $\mathcal{S}_A(x)$ and $\mathcal{S}_{A_\ell}(x)$ in Lemma 5.1 but also from the error term in (\ref{mainupper}).  Although the error term was negligable in $x$ with an appropriate choice of sieve level $D$, the method with which we establish the asymptotic in Theorem 2.3 will require subconvexity estimates for twisted Rankin-Selberg $L$-functions if one wishes to show that the error term is smaller than the main term with respect to the spectral parameter.}
\end{Remark2.3}

In \S4.1, Lemma 4.1 will show that by a ``partial Sato-Tate'' result one has
\begin{equation*}
M_u(x)\ll_u (\log z)^{-1/6}.
\end{equation*}
This fractional power saving of $\log z$ comes from the holomorphy and non-vanishing of the associated symmetric square, fourth and sixth power $L$-functions up to the line $\textnormal{Re}(s)=1$.  With the above bound for $M_u(x)$, we have shown
\begin{equation*}
\sum_{n \leqslant x} | \lambda ( n ) \lambda ( n + \ell ) | \ll_{\ell}\frac{x}{(\log z)^{1/6}} L(\textnormal{sym}^2 u, 1)\ll_{\ell,u} \frac{x}{(\log x)^{1/7}}
\end{equation*}
after the appropriate choice for $z=x^{1/s}$ to control the size of the error term $R_u(x)$. This concludes the proof of Theorem 1 with $\delta=1/7$.

\newpage
\section{Rankin-Selberg in arithmetic progression}
Let $\ell \neq 0$ be a small fixed integer and let $z=x^{1/s}$ for some $s$ to be chosen later. Let $a, a_\ell$ be small fixed integers
such that $( a, a_\ell )=(a a_\ell,\ell)=1$ and $p|aa_\ell \Rightarrow p \leqslant z$ and $d, d_\ell$
fixed integers such that $d, d_\ell |P_{6\ell} ( z )$ where
\begin{equation*} 
P_{6\ell} ( z ) := \prod_{\substack{
     p \leqslant z\\
     p \nmid 6\ell
   }} p. 
\end{equation*}
A key part in the proof of the main theorem will be to establish an asymptotic, when the linear equation $a_\ell d_\ell c_\ell = adc + \ell$ is solvable in $c$ and $c_\ell$, for the partial sums
\begin{equation*} 
 A_{d, d_\ell} ( x ) := \sum_{\substack{
     c \leqslant x/ad\\
     a_\ell d_\ell c_\ell= a d c + \ell
   }} \eta ( d c ) 
\end{equation*}
with $\eta$ a non-negative, completely multiplicative function defined by
\begin{equation*}
\eta(p^\alpha):=\left\{\begin{array}{ll}
\lambda^2(p) & \textnormal{if }  p \mid 6 \\
\lambda^{2\alpha}(p) & \textnormal{if }  p \nmid 6
\end{array}\right.
\end{equation*}
for prime powers $p^\alpha$. Here $\lambda ( p )$ is the $p$-th Hecke eigenvalue for some fixed Hecke-Maass cusp form $u$.  Observe that our choice of integers $a, a_\ell, d, d_\ell$ ensures that $(ad,a_\ell d_\ell)=(a d a_\ell d_\ell,\ell)=1$ when the linear equation is solvable.

Since we would like to use this asymptotic in conjunction with an upper-bound
sieve application, we want it to be of the form
\begin{equation*} 
A_{d, d_\ell} ( x ) = g' ( d ) g'' ( d_\ell ) X + r_{d, d_\ell} ( x ) 
\end{equation*}
where $g' ( d ) g'' ( d_\ell ) X$ is the expected main term with $g'$ and $g''$ multiplicative functions and the error term
$r_{d, d_\ell} ( x )$ is sufficiently small.  Here $X$ is approximately equal to 
\begin{equation*}
A(x)=\sum_{\substack{ab\leqslant x\\ a_\ell b_\ell = a b + \ell}} \eta(b).
\end{equation*}
In order to achieve this task, we interperet the linear equation $a_\ell d_\ell c_\ell=a d c + \ell$ as the congruence
\begin{equation*} c \equiv - \ell \overline{ad} ( \textnormal{mod } a_\ell d_\ell ) \end{equation*}
and employ some standard results from Dirichlet series analysis.
\begin{lemma3a}
  Let $D ( s )$ be a Dirichlet series with positive coefficients $a_n$ which
  converges absolutely for $\textnormal{Re} ( s ) > 1$ and extends to a meromorphic function in
  the half plane $\textnormal{Re} ( s ) \geqslant 1 / 2$ with a simple pole of
  residue $r$ at the point $s = 1$
  \begin{equation*} D ( s ) := \sum_{n \geqslant 1} \frac{a_n}{n^s} . \end{equation*}
  If D(s) satisfies the bound
 \begin{equation} D ( s ) \ll (|s|Q)^{1 + \varepsilon} \label{10.9}
  \end{equation}
  for some positive constant $Q$ and any $\varepsilon>0$ on the line $\textnormal{Re} ( s ) =
  1 / 2$, then for any $x \geqslant 1$ we have
  \begin{equation} \sum_{n \leqslant x} a_n = rx + O ( x^{3/4+\varepsilon}Q^{1/2+\varepsilon}) \label{10.10} \end{equation}
for any $\varepsilon>0$.
\begin{proof}
Choose a nice compactly supported test function $h(u)$ to smooth out our sum in question, which majorizes the dyadic segment $[x,2x]$ and is supported on $[x-y,2x+y]$ for some $y\leqslant x$ to be chosen later so that
\begin{equation*}
\sum_{n} a_{n} h (n)= \sum_{x-y \leqslant n < x}  a_{n} h (n) +  \sum_{x \leqslant n < 2x}  a_{n} +  \sum_{2x \leqslant n < 2x+y}  a_{n} h(n).
\end{equation*}   
Using the Mellin transform of $h(u)$, we start with 
\begin{equation*} 
   \sum_n a_n h ( n ) = \frac{1}{2 \pi i} \int_{( 2 )} H ( s ) D ( s ) ds. 
\end{equation*}
and move the line of integration to $\textnormal{Re} ( s ) = 1 / 2$ picking up the
pole at $s = 1$,
\begin{equation*}
\sum_n a_n h ( n ) = r H(1) + \frac{1}{2 \pi i} \int_{( \frac{1}{2} )} H ( s ) D( s ) ds. 
\end{equation*} 
Repeated partial integration and appropriate choice of $y$ then gives
\begin{equation*}
\sum_{n} a_n h ( n ) = rx + O(x^{3/4+\varepsilon}Q^{1/2+\varepsilon}).
\end{equation*}
If we had instead chosen a similar nice test function $g(u)$ supported on $[x,2x]$ with 
\begin{equation*}
\sum_{n} a_{n} g (n)\leqslant \sum_{x \leqslant n < 2x}  a_{n},
\end{equation*} 
then we would have
\begin{equation}
\sum_{x<n\leqslant 2x} a_n=rx+O(x^{3/4+\varepsilon}Q^{1/2+\varepsilon})\label{prooflemma3a}.
\end{equation}
The result follows after subdividing our original interval $[1,x]$ into dyadic segments and applying (\ref{prooflemma3a}) for each segment.
\end{proof}
\end{lemma3a}

In order to apply this lemma to our setting, we construct Dirichlet series
\begin{equation}
D_{\chi}(s):=\sum_{n \geqslant 1} \frac{\chi(n) \eta(n)}{n^{s}} \label{Zchi}
\end{equation}
for every character $\chi$ of a fixed modulos $q$ and note its properties.  We start first, however, with the full series
\begin{equation}
D_u(s):=\sum_{n\geqslant 1}  \frac{\eta(n)}{n^{s}}=\prod_{p\mid 6}\bigg(1+\frac{\lambda^2(p)}{p^s}\bigg)\prod_{p\nmid 6}\bigg(1-\frac{\lambda^2(p)}{p^s}\bigg)^{-1}
\label{RSflat}
\end{equation}
as each $D_{\chi}(s)$ will have similar properties. Note that 
\begin{equation*}
D_u(s)=\gamma_u(s) L(u\otimes u,s)
\end{equation*}
for some small correction factor $\gamma_u(s)$ whose product 
\begin{eqnarray}
\gamma_u(s)& = &\prod_{p\mid 6}\left(1 -\frac{p^{2s}+p^{s} - (2p^{2s}-p^{s}-1)\lambda^2(p)+(p^{2s}-p^{s})\lambda^4(p)}{p^{4s}(1+p^{-s})}\right)\nonumber\\
& \times &\prod_{p\nmid 6}\left(1 +\frac{2p^s\lambda^2(p)-p^s-1}{p^{3s}(1-\lambda^2(p)p^{-s})(1+p^{-s})} \right) \label{gammaeuler}
\end{eqnarray}
controls the difference in the second degree terms of the Euler products for $D_u(s)$ and the Rankin-Selberg $L$-function
\begin{equation*}
L(u\otimes u,s)=\sum_{n\geqslant 1}  \frac{\lambda^2(n)}{n^{s}}=\frac{\zeta(s)}{\zeta(2 s)}L(\textnormal{sym}^2 u,s).
\end{equation*}

By (\ref{localfactors}) we can write $\lambda^2(p)=\alpha^2_p+2+\alpha^{-2}_p$ and have by (\ref{KIMSARNAK}) that $|\alpha^2_p|\leqslant p^{7/32}$.  Therefore, with $0\leqslant \lambda^2(p) \leqslant p^{7/32} + 2 + p^{-7/32}$ for every prime $p$, one can check that we trivially have
\begin{equation*}
\frac{3}{5\pi^2}<\gamma_u:=\gamma_u(1)<15
\end{equation*}
 for any form $u$.  Furthermore, we know that $D_u(s)$ satisfies the same convexity bound as the Rankin-Selberg $L$-function for $\textnormal{Re}(s)=1/2$, namely
\begin{equation}
D_u(s)\ll_u(|s|)^{1+\varepsilon}
\label{RSconvexitybound}
\end{equation}
for any $\varepsilon>0$.

\begin{lemma3b}
Let $\chi_0$ be the principal character of modulos $q$. We then have
\begin{equation}
\sum_{n\leqslant x} \chi_0(n) \eta(n) = \frac{\theta}{\zeta(2)} L(\textnormal{sym}^2 u,1) x + O_u(x^{3/4+\varepsilon} q^{1/2+\varepsilon})
\end{equation}
for any $\varepsilon>0$ where
\begin{equation*}
\theta=\gamma_u \prod_{\substack{p\mid q \\ p\mid 6}} \left(1+\frac{\lambda^2(p)}{p}\right)^{-1}\prod_{\substack{p\mid q \\ p\nmid 6}} \left(1-\frac{\lambda^2(p)}{p}\right)
\end{equation*}
with constant $\gamma_u$ satisfying
\begin{equation*}
\frac{3}{5\pi^2} < \gamma_u < 15. 
\end{equation*}
\begin{proof}
Apply Lemma 3a for the Dirichlet series $D_{\chi_0}(s)$ with 
\begin{eqnarray*}
a_n=\chi_0(n) \eta(n) & \textnormal{and} & Q=q 
\end{eqnarray*}
noting that 
\begin{equation*}
\sum_{(n,q)=1} \frac{\eta(n)}{n^s} = D_u(s) \prod_{\substack{p\mid q\\p\mid 6}}\bigg(1+\frac{\lambda^2(p)}{p^s}\bigg)^{-1}\prod_{\substack{p\mid q\\p\nmid 6}}\bigg(1-\frac{\lambda^2(p)}{p^s}\bigg)
\end{equation*}
and that 
\begin{equation*}
D_{\chi_0}(s) \ll_u (|s|q)^{1+\varepsilon}
\end{equation*} 
for any $\varepsilon>0$ on the line $\textnormal{Re}(s)=1/2$.
\end{proof}
\end{lemma3b}
\noindent For every non-principal character $\chi$ of modulus $q$, we have that $D_{\chi}(s)$ is holomorphic for $\textnormal{Re}(s)\geqslant 1/2$.  Modifying the arguments slightly in the proofs of Lemma 3a and Lemma 3b, we obtain a bound for non-principal character twists. 
\begin{lemma3c}
Let $\chi$ be a non-principal character of modulos $q$. We then have
\begin{equation}
\sum_{n\leqslant x} \chi(n) \eta(n) \ll_u x^{3/4+\varepsilon} q^{1/2+\varepsilon}
\end{equation}
for any $\varepsilon>0$.
\end{lemma3c}
Using the Dirichlet character orthogonality relation 
\begin{eqnarray}
\frac{1}{\varphi ( q )}  \sum_{\chi ( q )} \bar{\chi} ( m ) \chi
  ( n )  =  \left\{\begin{array}{ll}
    1 & \textnormal{if } n \equiv m ( \textnormal{mod } q )\\
    0 & \textnormal{otherwise} 
  \end{array}\right. & \textnormal{when} & (m,q)=1
\label{ortho}
\end{eqnarray}
we combine Lemma 3b and Lemma 3c to have the following result.
\begin{lemma3d}
Fix integers $m$ and $q$ such that $(m,q)=1$. We then have
\begin{equation*}
 \sum_{\substack{
    n \leqslant x\\
    n \equiv m ( \textnormal{mod } q )
  }} \eta ( n )= \frac{\theta}{\zeta(2)} L(\textnormal{sym}^2 u,1) \frac{x}{\varphi(q)} + O_u\left(x^{3/4+\varepsilon} q^{1/2+\varepsilon}\right) \label{lemma3dresult}
\end{equation*}
for any $\varepsilon>0$ where
\begin{equation*}
\theta=\gamma_u \prod_{\substack{p\mid q \\ p\mid 6}} \left(1+\frac{\lambda^2(p)}{p}\right)^{-1}\prod_{\substack{p\mid q \\ p\nmid 6}} \left(1-\frac{\lambda^2(p)}{p}\right)
\end{equation*}
with constant $\gamma_u$ satisfying
\begin{equation*}
\frac{3}{5\pi^2} < \gamma_u < 15. 
\end{equation*}
\end{lemma3d}
\begin{Remark2.1}
\textnormal{The error term in Lemma 3d depends on the spectral parameter. If one looks at the gamma factors for the twisted Rankin-Selberg $L$-function, then one sees that the true convexity bound for $D_u(s)$ on the line $\textnormal{Re}(s)=1/2$ depends on the spectral parameter $t_u$,}
\begin{equation*}
D_u(s) \ll q^{1+\varepsilon}\left(|s|^{1+\varepsilon}+(|s|t_u)^{1/2+\varepsilon}\right).
\end{equation*}
\textnormal{Results similar to Lemma 3a would then at best produce the modified error term}
\begin{equation*}
\sum_{\substack{
    n \leqslant x\\
    n \equiv m ( \textnormal{mod } q )
  }} \eta ( n )= \frac{\theta}{\zeta(2)} L(\textnormal{sym}^2 u,1) \frac{x}{\varphi(q)} + O\left(x^{1/2} q^{1+\varepsilon} t_u^{1/2+\varepsilon}\right)
\end{equation*}
\textnormal{for any $\varepsilon>0$.  If $x$ is of size $t_u$, we simply do not have the error term smaller than the main term unless we first assume a subconvexity bound for $D_u(s)$ for $\textnormal{Re}(s)=1/2$.}
\end{Remark2.1}

Returning to our sifting sums in arithmetic progressions $A_{d,d_\ell}(x)$, if we interperet the linear equation $a_\ell d_\ell c_\ell=a d c + \ell$ with $(ad,a_\ell d_\ell)=(a d a_\ell d_\ell,\ell)=1$ as the congruence
\begin{equation*} c \equiv - \ell \overline{ad} ( \textnormal{mod } a_\ell d_\ell ) \end{equation*}
then 
\begin{equation*}
  \sum_{\substack{
    adc \leqslant x\\
    a_\ell d_\ell c_\ell= adc + \ell
  }} \eta ( dc ) = \lambda^2(d) \sum_{\substack{
    adc \leqslant x\\
    c \equiv - \ell \overline{ad} ( \textnormal{mod } a_\ell d_\ell )
  }} \eta ( c )
\end{equation*}
since $\eta(c d)=\lambda^2(d)\eta(c)$ for $d\mid P_{6\ell}(z)$.  Therefore, by Lemma 3d we have
\begin{eqnarray*}
 \sum_{\substack{
    adc \leqslant x\\
    a_\ell d_\ell c_\ell= adc + \ell
  }} \eta ( dc )&=&\frac{\theta}{\zeta(2)} L(\textnormal{sym}^2 u,1) \frac{\lambda^2(d)}{\varphi(a_\ell d_\ell)} \frac{x}{a d}\\ 
& &+ O_u\left(\lambda^2(d)\left(\frac{x}{ad}\right)^{3/4+\varepsilon}(a_\ell d_\ell)^{1/2+\varepsilon}\right)
\end{eqnarray*}
for any $\varepsilon>0$ where $\theta$ is now
\begin{equation*}
\theta=\gamma_u \prod_{\substack{p\mid a_\ell \\ p\mid 6}} \left(1+\frac{\lambda^2(p)}{p}\right)^{-1}\prod_{\substack{p\mid a_\ell d_\ell \\ p\nmid 6}} \left(1-\frac{\lambda^2(p)}{p}\right)
\end{equation*}
with
\begin{equation*}
\frac{3}{5\pi^2} < \gamma_u < 15. 
\end{equation*}

Since $|\lambda^2(d) d^{-3/4}|\ll 1$ by the Kim-Sarnak bound (\ref{KIMSARNAK}), if $a,a_\ell,d,d_\ell\leqslant x^{1/16}$ we see that the error term is $O_u(x^{13/16+\varepsilon})$.
This brings us to an asymptotic for the partial sums $A_{d,d_\ell}(x)$ and proves Theorem 2.3.

\newpage
\section{Some bounds involving Hecke eigenvalues}
Let $u$ be a Hecke-Maass cusp form on $X={\rm SL}\sb 2(\mathbbm{Z}) \backslash \mathbbm{H}$. We associate with $u$ the $L$-function
\begin{equation} L ( u, s ) = \sum_{n \geqslant 1} \lambda ( n ) n^{- s} = \prod_p ( 1 -
   \lambda ( p ) p^{- s} + p^{- 2 s} )^{- 1} \label{Lfunc} \end{equation}
and write each local factor as
\begin{equation} 1 - \lambda ( p ) p^{- s} + p^{- 2 s} = ( 1 - \alpha_p p^{- s} ) ( 1 -
   \beta_p p^{- s} )  \end{equation}
so that
\begin{eqnarray}
  \alpha_p + \beta_p = \lambda ( p ) & \textnormal{and} & \alpha_p \beta_p = 1. 
\label{localfactors}
\end{eqnarray}
Since the Hecke operators are self-adjoint, we have that $\lambda ( p ) \in \mathbbm{R}$ and that either $| \alpha_p | = |
\beta_p | = 1$ or $\alpha_p, \beta_p \in \mathbbm{R}$.  The current best
known estimate for these local factors is due to Kim and Sarnak [K-Sa]
\begin{equation} | \alpha_p |, | \beta_p | \leqslant p^{7 / 64} . \label{KIMSARNAK} \end{equation}
The Ramanujan-Petersson Conjecture would give $|\lambda(p)|\leqslant 2$.

In addition to $L ( u, s )$, we have the symmetric power $L$-functions
\begin{equation} L ( \textnormal{sym}^m u, s ) := \prod_p \prod^m_{j = 0} ( 1 - \alpha^{m -
   j}_p \beta^j_p p^{- s} )^{- 1} \label{Lsymm} \end{equation}
for $m\geqslant 1$ and the Rankin-Selberg convolution
\begin{equation} L(u  \times u) := \zeta (2s) L ( u \otimes u, s ) := \zeta ( 2 s ) \sum_{n \geqslant 1} \lambda^2 ( n )
   n^{- s} = \zeta(s) L(\textnormal{sym}^2 u, s) \label{RSLfunc} \end{equation}
which has a simple pole at $s = 1$. By the
works of Kim-Shahidi [K-Sh] and Kim [Ki], the symmetric $m$-th power $L$-functions for $m \leqslant 8$ are holomorphic and non-vanishing for  $\textnormal{Re}(s)\geqslant 1$.

We already know from Rankin-Selberg Theory that 
\begin{equation}
\sum_{n\leqslant x} \lambda^2(n) \sim c_u x
\label{RS}
\end{equation}
for some constant $c_u$ depending on our cusp form $u$.  Results involving the average value of $|\lambda(n)|$ were established in the works of several authors including Rankin ([Ra1], [Ra2]), Elliott, Moreno, Shahidi [E-M-S] and Murty [Mu].  Their works were based on properties of the symmetric square and symmetric fourth power $L$-functions associated with $u$ at the point $s=1$.  In particular, the work of Elliott, Moreno and Shahidi [E-M-S] established the bound
\begin{equation}
\sum_{n\leqslant x} |\lambda(n)| \ll_u \frac{x}{(\log x)^{1/18}}
\label{EMS}
\end{equation}
under the assumption of the Ramanujan-Petersson conjecture.  The above bound follows from the simple inequality
\begin{eqnarray}
 | y | & \leqslant & 1 + \frac{1}{2} ( y^2 - 1 ) - \frac{1}{18} ( y^2 - 1 )^2 \nonumber \\ 
& = & \frac{17}{18} + \frac{11}{18} ( y^2 - 1 ) - \frac{1}{18} ( y^4 - 2 )\label{RPineq},
\end{eqnarray}
which holds for all real $y$ with $|y|\leqslant2$, and from the Hecke relations
\begin{eqnarray}
\lambda^2(p)-1 & = & (\alpha^2_p+1+\beta^2_p) = \lambda(p^2), \label{sym2eq}\\
\lambda^4(p)-2 & = & (\alpha^4_p+\alpha^2_p+1+\beta^2_p+\beta^4_p)+3(\alpha^2_p+1+\beta^2_p)\label{sym4eq}\\
& = & \lambda(p^4)+3\lambda(p^2). \nonumber
\end{eqnarray}
We refer the reader to that paper for details and instead apply their ideas to properly bound $M_u(x)$ in Lemma 1.3.

\subsection{Saving a fractional power of  $\log z$ from $M_u(x)$}
We would like to establish a bound of the form 
\begin{equation}
M_u(x)\ll_u (\log z)^{-\delta}
\end{equation}
for some $\delta>0$ without the assumption of the Ramanujan-Petersson conjecture.  This will require us to construct an inequality similar to (\ref{RPineq}), but which holds for all real values $y$.  To this end, we introduce an extra term and look for values $a$ and $b$ such that
\begin{equation}
|y| \leqslant 1+\frac{1}{2}(y^2-1)+a(y^2-1)^2+b(y^2-1)^3 \label{MYineqform}
\end{equation}
holds true for all real values of $y$, for then we would have
\begin{equation*}
-(1-|y|)^2=2|y|-y^2-1 \leqslant  2a (y^2-1)^2 + 2b (y^2-1)^3.
\end{equation*}
One such possible choice is $(a,b)=(-1/9, 1/36)$ giving
\begin{equation*}
-(1-|y|)^2 \leqslant -\frac{3}{18}+\frac{11}{18}(y^2-1)-\frac{7}{18}(y^4-2)+\frac{1}{18}(y^6-5).
\end{equation*}
Note that the right hand side of the above inequality is strictly less than 0 for all $|y|\leqslant 2$ except for $|y|=1$.

Substituting $y=\lambda(p)$ for each prime $p$ and averaging over primes up to $z$ we get
\begin{equation}
- \sum_{p\leqslant z} \frac{(1-|\lambda(p)|)^2}{p} \leqslant -\frac{3}{18} \log \log z + O_u(1) \label{absmoment5}
\end{equation}
provided that 
\begin{equation}
\frac{11}{18}\sum_{p\leqslant z}\frac{\lambda^2(p)-1}{p} - \frac{7}{18}\sum_{p\leqslant z}\frac{\lambda^4(p)-2}{p}+\frac{1}{18}\sum_{p\leqslant z}\frac{\lambda^6(p)-5}{p}\label{symmpartial}
\end{equation}
can be bounded independently of $z$.  From the Hecke relations (\ref{sym2eq}), (\ref{sym4eq}) and 
\begin{eqnarray}
\lambda^6(p)-5 & = & (\alpha^6_p+\alpha^4_p+\alpha^2_p+1+\beta^2_p+\beta^4_p+\beta^6_p) \nonumber\\
& & + 5(\alpha^4_p+\alpha^2_p+1+\beta^2_p+\beta^4_p)+9(\alpha^2_p+1+\beta^2_p),\label{sym6eq}\nonumber\\
& = & \lambda(p^6) + 5 \lambda(p^4) +9 \lambda(p^2)
\end{eqnarray}
one sees that (\ref{symmpartial}) may be rewritten as
\begin{equation*}
-\frac{1}{18}\sum_{p\leqslant z}\frac{\lambda(p^2)}{p} - \frac{2}{18}\sum_{p\leqslant z}\frac{\lambda(p^4)}{p}+\frac{1}{18}\sum_{p\leqslant z}\frac{\lambda(p^6)}{p}
\end{equation*}
and is related to the logarithms of the associated symmetric square, fourth and sixth power $L$-functions at the point $s=1$ as is evident by the Euler product definition (\ref{Lsymm}).  Since the holomorphicity and non-vanishing of these symmetric powers is known for $\textnormal{Re}(s)\geqslant 1$, we have (\ref{absmoment5}) and can state the following Lemma.
\begin{lemma4.1}
For any Hecke-Maass cusp form $u$ with Hecke eigenvalues $\lambda$, we have for $M_u(x)$ as given in (\ref{MTHM3})
\begin{equation}
M_u(x)\ll \left(\frac{L_6(u,z)}{L_2(u,z) L^2_4(u,z) (\log z)^3}\right)^{1/18}\ll_u (\log z)^{-1/6} \label{lemma4.1bresult}
\end{equation}
where the
\begin{equation}
L_{m}(u,z)  :=  \prod_{p\leqslant z} \prod_{j=0}^m\left(1-\alpha_p^{m-j}\beta_p^{j} p^{-1}\right)^{-1}\label{symmz}
\end{equation}
are the partial Euler products of the associated symmetric power $L$-functions at the point $s=1$.
\end{lemma4.1}

\subsection{Fourth moment bound in terms of symmetric powers}
By the Kim-Sarnak bound (\ref{KIMSARNAK}) we have that
\begin{equation*}
\sum_{n \leqslant x} \lambda^4 ( n ) \ll x \prod_{p\leqslant x}\left(1+\frac{\lambda^4(p)}{p}\right)\ll x \exp\left(\sum_{p\leqslant x}\frac{\lambda^4(p)-2}{p}\right)\log^2(x).
\end{equation*}
Applying the Hecke relations and our partial product definitions $L_m(u,z)$ in (\ref{symmz}) we get the following bound.
\begin{lemma4.2}
  Let $u$ be a Hecke-Maass cusp form with Hecke eigenvalues $\lambda$, then
  \begin{equation} 
\sum_{n \leqslant x} \lambda^4 ( n ) \ll x ( \log x )^2 L_4(u,x) L^3_2(u,x). 
\label{lemma4.2bound}
\end{equation}
\end{lemma4.2}
\noindent This bound will be used in \S 5.1 in the proof of Lemma 5.1.

\newpage
\section{Technical Lemmata}
During our proof of Theorem 1 in \S 2, we stated that the sums $\mathcal{S}_A(x)$ and $\mathcal{S}_{A_\ell}(x)$ defined in (\ref{SA}) and (\ref{SAl}) would have small contribution to our shifted sum $S_{\ell}(x)$. We also removed square-full integers $b$ from our sifting sequence in \S 2.3 in order to define a new multiplicative function $\eta$ and only applied the double upper-bound sieve to the sum $S_{a,a_\ell}(x,z)$ defined in (\ref{sievingsum}). We now show the necessary Lemmata which justify our analysis.
\subsection{$\mathcal{S}_A(x)$ and $\mathcal{S}_{A_\ell}(x)$ have small contribution}
With our factorization of $n$ and $n+\ell$ as
\begin{eqnarray*}
  n = ab & \textnormal{and} & n + \ell = a_\ell b_\ell,
\end{eqnarray*}
such that for every prime $p$ dividing $n(n+\ell)$,
\begin{eqnarray*}
 p| a a_\ell  \Rightarrow p \leqslant z & \textnormal{and} & (b b_\ell, P(z))=1,
\end{eqnarray*}
we partitioned the shifted sum $S_{\ell}(x)$ based on the size of $a$ and $a_\ell$. We defined
\begin{eqnarray*}
  \mathcal{S}_A ( x ) & = & \sum_{\substack{
    n = ab \leqslant x\\
    n + \ell = a_\ell b_\ell\\
    p| ( aa_\ell ) \Rightarrow p \leqslant z\\
    (b b_\ell, P(z))=1\\
    a > x^{1 / 16}
  }} | \lambda ( n ) \lambda ( n + \ell ) |,\\
  \mathcal{S}_{A_\ell} ( x ) & = & \sum_{\substack{
    n = ab \leqslant x\\
    n + \ell = a_\ell b_\ell\\
    p| ( aa_\ell ) \Rightarrow p \leqslant z\\
    (b b_\ell, P(z))=1\\
    a_\ell > x^{1 / 16}
  }} | \lambda ( n ) \lambda ( n + \ell ) |,
\end{eqnarray*}
and now prove the following Lemma using a trick of Rankin to detect numbers with small prime factors.
\begin{lemma5.1}
For a fixed Hecke-Maass cusp form $u$ and $\ell\neq 0$ we have
\begin{equation*}
\mathcal{S}_A(x)+\mathcal{S}_{A_\ell}(x) \ll_{u,\ell} x (\log x)^{-C}
\end{equation*}
for any constant $C > 0$.
\end{lemma5.1}

Define
\begin{equation*} 
\Phi ( x, z ) := \sum_{\substack{
     a \leqslant x\\
     p|a \Rightarrow p \leqslant z }} 1 
\end{equation*}
to be the number of integers less than $x$ with prime factors less than $z$. For any $\alpha>0$, we have
\begin{equation*}
\Phi (x,z) \leqslant \sum_{\substack{a \leqslant x \\ p|a \Rightarrow p \leqslant z }}
\left(\frac{x}{a} \right)^{\alpha} \leqslant x^{\alpha} \prod_{p\leqslant z} \left( 1-\frac{1}{p^{\alpha}}\right)^{-1}.
\end{equation*}
Setting $\alpha=1-\eta$ with $\eta \longrightarrow 0$ as $z \longrightarrow \infty$ gives
\begin{eqnarray*}
\Phi (x,z) & \ll & x^{\alpha} \prod_{p\leqslant z} \left( 1+\frac{1}{p^{\alpha}}\right)\\
 & \ll & x \exp(-\eta \log x + \sum_{p\leqslant z} \frac{1}{p} (1+(\eta \log p) z^{\eta}))\\
 & \ll & x (\log z) \exp (-\eta \log x +\eta z^\eta (\log z)).
\end{eqnarray*}
Choosing $\eta = (\log z)^{-1}$ we get the bound
\begin{equation}
\Phi (x,z) \ll x (\log z) \exp(-\frac{\log x}{\log z}). \label{Rank1}
\end{equation}
Better bounds may be established (see [deBr]), however, (\ref{Rank1}) is sufficient for our application.

In addition to (\ref{Rank1}) we have
\begin{equation} \Psi ( x, z ) := \sum_{\substack{
     b \leqslant x\\
     (b,P(z))=1
   }} 1 \ll x \prod_{p \leqslant z} \left( 1 - \frac{1}{p} \right) \ll
   \frac{x}{\log z} . \label{Psi} \end{equation}
Applying (\ref{Rank1}), (\ref{Psi}) and Lemma 4.2 to $\mathcal{S}_A( x )$ allows us to treat this sum immediately. 

By H\"{o}lder's inequality and Lemma 4.2 we get
\begin{eqnarray*}
  \mathcal{S}_A ( x ) & \leqslant & \bigg( \sum_{\substack{
    x^{1 / 16} < a \leqslant x\\
    p|a \Rightarrow p \leqslant z\\
  }} \sum_{\substack{
    b \leqslant x / a\\
    (b, P(z))=1
  }} 1 \bigg)^{1 / 2} \bigg( \sum_{n \leqslant x} | \lambda ( n ) |^4 \bigg)^{1 /
  4} \bigg( \sum_{n \leqslant x} | \lambda ( n + \ell ) |^4 \bigg)^{1 / 4}\\
  & \ll_{\ell} & \bigg( \frac{x}{\log z}  \sum_{\substack{
    x^{1 / 16} < a \leqslant x\\
    p|a \Rightarrow p \leqslant z
  }}  \frac{1}{a} \bigg)^{1 / 2}  \bigg( x ( \log x )^2 L_4(u,x) L_2^3(u,x)\bigg)^{1 / 2} 
\end{eqnarray*}
with $L_2(u,x)$ and $L_4(u,x)$ as defined in (\ref{symmz}).  After partial summation,
\begin{equation*}
  \nonumber \sum_{\substack{
    x^{1 / 16} < a \leqslant x\\
    p|a \Rightarrow p \leqslant z
  }}  \frac{1}{a} 
\ll  ( \log z ) ( \log x ) \exp ( - \frac{1}{16}  \frac{\log x}{\log z} )
\end{equation*}
and combining our results establishes the bound
\begin{equation*}
\mathcal{S}_A(x) \ll_{\ell} x  (\log x)^{3/2} \exp(-\frac{1}{32} \frac{\log x}{\log z})L_4^{1/2}(u, x) L_2^{3/2}(u,x).
\end{equation*}
After taking $z = x^{1 / s}$ we get
\begin{equation*} 
\mathcal{S}_A ( x ) \ll_{u,\ell} x ( \log x )^{-C} 
 \end{equation*}
for any $C>0$ by choosing $s=c \log \log x$ appropriately. Applying the same analysis to $\mathcal{S}_{A_\ell}(x)$ concludes the proof of Lemma 5.1.

\subsection{Square-full $b$ have small contribution}
In this section, we justify the restriction to the sum over square-free $b$ in the proof of Theorem 1.  Recall that in \S 2.3 we first removed square-full integers $b$ and then took advantage of positivity by setting our sifting sequence to be $\mathcal{A}=(\eta(b))$ for some new multiplicative function $\eta$ defined in (\ref{eta}) claiming one only need apply the double upper-bound sieve to the sum (\ref{sievingsum}) as the contribution from square-full $b$ would be small.

Starting with the sum over all $b$ for fixed $a,a_\ell\leqslant x^{1/16}$ with $(a,a_\ell)=1$, we write 
\begin{equation*}
\sum_{\substack{
     b \leqslant x / a\\
     a_\ell b_\ell = ab + \ell\\
     ( bb_\ell, P ( z ) ) = 1
   }} \lambda^2 ( b ) = \sum_{\substack{
     b \leqslant x / a\\
     a_\ell b_\ell = ab + \ell\\
     ( bb_\ell, P ( z ) ) = 1
   }} \lambda^2 ( b ) \mu^2(b)  + \mathfrak{S}_{a,a_\ell}(x,z)
\end{equation*}
where $\mathfrak{S}_{a,a_\ell}(x,z)$ is the sum over those $b$ which contain square factors.  If we simply pull out the factor $m$ from $b$ which consists of all prime power factors $p^\alpha$ where $\alpha\geqslant 2$ we have 
\begin{equation*}
\mathfrak{S}_{a,a_\ell}(x,z)=\sum^{\sharp\sharp}_{\substack{
     z^2 \leqslant m \leqslant x / a\\
     (m,P(z))=1}}
     \lambda^2 ( m )
\sum_{\substack{
     b \leqslant x / am\\
     a_\ell b_\ell = amb + \ell\\
     (b,m)=1\\
     ( bb_\ell, P ( z ) ) = 1
   }} \lambda^2 ( b ) \mu^2(b).
\end{equation*}
The $\sharp \sharp$ super-script means that each prime power factor $p^\alpha$ dividing $m$ has $\alpha\geqslant 2$.  

Dropping the co-primality conditions for $(b b_\ell,P(z))=1$ and $(b,m)=1$ by positivity while maintaining that $(a_\ell, am)=1$, we can apply Lemma 3d, since $\lambda^2(b)\mu^2(b) \leqslant \eta(b)$ for all $b$, to have
\begin{equation*}
\sum_{\substack{
     b \leqslant x / am\\
     a_\ell b_\ell = amb + \ell
   }} \lambda^2 ( b ) \mu^2(b) \leqslant \frac{\theta}{\zeta(2)} L(\textnormal{sym}^2 u,1) \frac{x}{\varphi(a_\ell) am} + O_u\left(\left(\frac{x}{am}\right)^{3/4+\varepsilon} a_\ell^{1/2+\varepsilon}\right)
\end{equation*}
for any $\varepsilon>0$ where $\theta$ is now
\begin{equation*}
\theta=\gamma_u \prod_{\substack{p\mid a_\ell \\ p\mid 6}} \left(1+\frac{\lambda^2(p)}{p}\right)^{-1}\prod_{\substack{p\mid a_\ell \\ p\nmid 6}} \left(1-\frac{\lambda^2(p)}{p}\right)
\end{equation*}
with
\begin{equation*}
\frac{3}{5\pi^2} < \gamma_u < 15. 
\end{equation*}
Noting that the Kim-Sarnak bound (\ref{KIMSARNAK}) and partial summation gives
\begin{equation*}
\sum^{\sharp\sharp}_{\substack{
     z^2 \leqslant m \leqslant x / a\\
     (m,P(z))=1}}
     \frac{\lambda^2 ( m )}{m^{3/4}} \ll z^{-1/32},
\end{equation*}
we see that
\begin{equation}
\mathfrak{S}_{a,a_\ell}(x,z) \ll_u \frac{x}{a a_\ell z^{1/32}}.\label{mathfrakbound}
\end{equation}
We therefore have the following Lemma.
\begin{lemma5.2}
For $\mathfrak{S}_{a,a_\ell}(x,z)$ as defined above, we have
\begin{equation*}
\sum_{\substack{
     a, a_\ell \leqslant x^{1 / 16}\\
     p|aa_\ell \Rightarrow p\leqslant z\\
     (a,a_\ell)=1
   }} \frac{| \lambda ( a ) \lambda ( a_\ell ) |}{ a a_\ell}\mathfrak{S}_{a,a_\ell}(x,z) \ll_u \frac{x}{z^{1/64}}
\end{equation*}
\begin{proof}
We insert the bound (\ref{mathfrakbound}) to get
\begin{equation*}
\sum_{\substack{
     a, a_\ell \leqslant x^{1 / 16}\\
     p|aa_\ell \Rightarrow p\leqslant z\\
     (a,a_\ell)=1
   }} \frac{| \lambda ( a ) \lambda ( a_\ell ) |}{ a a_\ell}\mathfrak{S}_{a,a_\ell}(x,z) \ll_u \frac{x}{z^{1/32}} \prod_{p\leqslant z} \left(1+\frac{|\lambda(p)|}{p}\right)^2.
\end{equation*}
By arguing as in \S4.1 we get the desired result.
\end{proof}
\end{lemma5.2}
\noindent Therefore the contribution from square-full $b$ is indeed small.

\newpage
\appendix
\section{Reducing to shifted convolution sums}
We demonstrate how one may relate the inner products $<f u_j, u_j>$ in the question of QUE to the study of a controlled number(independent of $t_j$) of shifted convolution sums of Hecke eigenvalues.
\subsection{Fourier coefficients}
For $f$ a fixed cusp form or incomplete Eisenstein series, we can express $f$
as a Fourier series expansion of the type
\begin{equation} f ( z ) = a_0 ( y ) + \sum_{\ell \neq 0} a_{\ell} ( y ) e ( \ell x ), \label{2.3} \end{equation}
with $a_0 ( y ) = 0$ in the cusp form case.  If $f ( z )$ is a fixed cusp
form with $\Delta$ eigenvalue $1/4+r^{2}$, then we have the expansion
\begin{equation*} f ( z ) = \sqrt{y}  \sum_{\ell \neq 0} \rho ( \ell ) K_{ir} ( 2 \pi | \ell
   |y ) e ( \ell x ) \end{equation*}
where the $\rho ( \ell )$ are complex numbers.  For $f ( z ) = E ( z| \psi )$
the incomplete Eisenstein series
\begin{equation*} E ( z| \psi ) : = \sum_{\gamma \in \Gamma_{\infty} \backslash \Gamma} \psi (
   \textnormal{Im } \gamma z ), \end{equation*}
where $\psi ( y )$ is a smooth function, compactly supported on
$\mathbbm{R}^+$, the coefficients $a_{\ell} ( y )$ can be determined in terms
of the coefficients of the Eisenstein series
\begin{equation*} E ( z, s ) : = \sum_{\gamma \in \Gamma_{\infty} \backslash \Gamma} ( \textnormal{Im }
   \gamma z )^s . \end{equation*}
The latter has the Fourier expansion
\begin{equation} E ( z, s ) = y^s + \varphi ( s ) y^{1 - s} + \sqrt{y}  \sum_{\ell \neq 0}
   \varphi_{| \ell |} ( s ) K_{s - \frac{1}{2}} ( 2 \pi | \ell |y ) e ( \ell x
   ) \label{2.4} \end{equation}
where
\begin{eqnarray*}
  \varphi ( s ) = & \sqrt{\pi}  \frac{\Gamma ( s - \frac{1}{2} ) \zeta ( 2 s -
  1 )}{\Gamma ( s ) \zeta ( 2 s )} = & \frac{\theta ( 1 - s )}{\theta ( s
  )},\\
  \theta ( s ) = & \pi^{- s} \Gamma ( s ) \zeta ( 2 s ), & \\
  \varphi_{\ell} ( s ) = & \frac{2}{\theta ( s )}  \sum_{ab = \ell} \left(
  \frac{a}{b} \right)^{s - \frac{1}{2}}, & \textnormal{if } \ell \geqslant 1.
\end{eqnarray*}
Note that
\begin{equation}
  \textnormal{res}_{s = 1} E ( z, s ) = \textnormal{res}_{s = 1} \varphi ( s ) = \frac{3}{\pi} . \label{2.5}
\end{equation}
We have
\begin{equation} E ( z| \psi ) = \frac{1}{2 \pi i} \int_{( 2 )} \Psi ( - s ) E ( z, s ) ds 
   \label{2.5.5} \end{equation}
where $\Psi ( s )$ is the Mellin transform of $\psi ( y )$.  This is an entire
function with rapid decay in vertical strips, specifically
\begin{equation*} \Psi ( s ) \ll ( |s| + 1 )^{- A} \end{equation*}
for any $A \geqslant 0$, $- 2 \leqslant \textnormal{Re} ( s ) \leqslant 2$, with the
implied constant depending only on $\psi$ and $A$.

From these formulas we find the coefficients of (\ref{2.3}) in the incomplete
Eisenstein series case,
\begin{equation*} a_0 ( y ) = \frac{1}{2 \pi i} \int_{( 2 )} \Psi ( - s ) ( y^s + \varphi ( s
   ) y^{1 - s} ) ds = \psi ( y ) + O ( y^{- 1} ) \end{equation*}
and for $\ell \neq 0$ we move the integration to the line $\textnormal{Re}(s) = 1/2$ to get
\begin{equation} a_{\ell} ( y ) = \left( \frac{y}{\pi} \right)^{\frac{1}{2}} \int^{+
   \infty}_{- \infty} \frac{\pi^{it} \Psi ( - \frac{1}{2} - it )}{\Gamma (
   \frac{1}{2} + it ) \zeta ( 1 + 2 it )}  \left( \sum_{ab = \ell} \left(
   \frac{a}{b} \right)^{it} \right) K_{it} ( 2 \pi | \ell |y ) dt. \label{2.6}
\end{equation}
Doing the same for $a_0 ( y )$ we get
\begin{equation*} a_0 ( y ) = \frac{3}{\pi} \Psi ( - 1 ) + O ( \sqrt{y} ). \end{equation*}
On the other hand, by unfolding the incomplete Eisenstein series $f$, we
derive that
\begin{equation*} < f, 1 > = \int^{1 / 2}_{- 1 / 2} \int_0^{\infty} \psi ( y ) d \mu z = \Psi
   ( - 1 ). \end{equation*}
Therefore, we have
\begin{equation} a_0 ( y ) = \frac{3}{\pi} < f, 1 > + O ( \sqrt{y} )  \label{zerocoeffEisenstein}\end{equation}
which conveniently contains the expected main term $3/\pi <f,1>$ in (\ref{1}).
For $\ell \neq 0$ we employ (\ref{B.1}) in Appendix B to (\ref{2.6}) in order to obtain
\begin{equation*} a_{\ell} ( y ) \ll \tau ( | \ell | ) | \ell |^{- 2} y^{- \frac{3}{2}} . \end{equation*}
Here $\tau ( \ell )$ is the divisor function. In the case of the cusp form, we
get a similar bound
\begin{equation} a_{\ell} ( y ) = \rho ( \ell ) \sqrt{y} K_{ir} ( 2 \pi | \ell |y ) \ll |
   \rho ( \ell ) || \ell |^{- 2} y^{- \frac{3}{2}} ( r^2 + 1 ) . \label{2.8a} \end{equation}
We state the results in the following Lemma.
\begin{lemmaA.1}
  Let $f \in \mathcal{A} ( X )$ be an automorphic function on $X =
 {\rm SL}\sb 2(\mathbbm{Z}) \backslash \mathbbm{H}$ with Fourier series expansion
  \begin{equation*} f ( z ) = a_0 ( y ) + \sum_{\ell \neq 0} a_{\ell} ( y ) e ( \ell x ) . \end{equation*}
  If $f$ is a cusp form in $\mathcal{A}_s ( X )$, then
  \begin{eqnarray*}
    & & a_0 ( y ) = 0  \textnormal{ and}\\
    & & a_{\ell} ( y ) \ll | \rho ( \ell ) || \ell |^{- 2} y^{- \frac{3}{2}} (
    |s|^2 + 1 )  \textnormal{ for }  \ell \neq 0.
  \end{eqnarray*}
 For $f$ an incomplete Eisenstein series,
  \begin{eqnarray*}
    & & a_0 ( y )  =  \frac{3}{\pi} < f, 1 > + O ( \sqrt{y} ) \textnormal{ and}\\
    & & a_{\ell} ( y ) \ll  \tau ( | \ell | ) | \ell |^{- 2} y^{- \frac{3}{2}}  \textnormal{ for }  \ell \neq 0.
  \end{eqnarray*}
\end{lemmaA.1}

One should note, that the bound (\ref{2.8a}) depends heavily on the spectral
parameter $r$.  When $f$ is a fixed cusp form during analysis, this will cause
no problem in estimations as $r$ will be a fixed constant.  Much more careful
attention has to be paid when we will be appealing to the Fourier series
expansions of the variable cusp forms $u_j$. Estimates must here be explicit and quite precise in terms of the spectral
parameter $t_j$.  

By ([Iw], 8.7) we have
\begin{equation} \sum_{|n| \leqslant x} | \rho_j ( n ) |^2 \ll ( t_j + x ) e^{\pi t_j} \label{2.9a} \end{equation}
where the implied constant is absolute.  The second moment (\ref{2.9a}), will make
an appearance when our analysis forces us to finally consider the inner
product $< u_j, u_j >$ in terms of the Fourier expansions for the variable
cusp forms $u_j$. At that point, the introduction of the spectral parameter is unavoidable.

\subsection{Transformations by unfolding method}

Let $u ( z )$ be a cusp form with $\Delta$ eigenvalue $\lambda = 1 / 4 + r^2$ and $f (
z )$ an automorphic function which is smooth and bounded on $\mathbbm{H}$. We
are going to compute the inner product $< fu, u >$ asymptotically in two ways
by an unfolding technique.  To this end, we fix a function $g ( y ) \in
C^{\infty}_c (\mathbbm{R}^+ )$, smooth and compactly supported on
$\mathbbm{R}^+$, and let
\begin{equation*} G ( s ) := \int^{+ \infty}_0 g ( y ) y^{s - 1} dy \end{equation*}
be its Mellin transform.  Therefore, $G ( s )$ is entire and
\begin{equation*} G ( s ) \ll ( |s| + 1 )^{- A} \end{equation*}
for any $A \geqslant 0$, uniformly in vertical strips, where the implied
constant depends on $g$ and $A$.

Let $X \geqslant 1$ and consider the integral
\begin{equation}
 I_f ( X ) := \frac{1}{2 \pi i} \int_{( \sigma )} G ( - s ) X^s
   \int_{\Gamma \backslash \mathbbm{H}} E ( z, s ) f ( z ) |u ( z ) |^2 d \mu zds \label{3.2} \end{equation}
with $\sigma > 1$.   Unfolding the inner integral we get
\begin{eqnarray}
  \nonumber I_f ( X ) & = & \frac{1}{2 \pi i} \int_{( \sigma )}  G ( - s ) X^s \int_{0}^{\infty}\int_{-1/2}^{1/2}y^s f ( z
  ) |u ( z ) |^2 d \mu zds\\
  & = & \int^{\infty}_0 g ( Xy ) y^{- 2} \left( \int^{1 / 2}_{- 1 / 2} f ( z
  ) |u ( z ) |^2 dx \right) dy. \label{3.3}
\end{eqnarray}
We show that our integral $I_f ( X )$ is bounded by $X$ by using the
following Lemma([Iw], Lemma 2.10).
\begin{lemmaIwaniec}
  Let $z \in \mathbbm{H}$ and $Y > 0$. We have
  \begin{equation*} \# \left\{ \gamma \in \Gamma_{\infty} \backslash \Gamma \, | \, \textnormal{Im} (
     \gamma z ) > Y \right\} < 1 + \frac{10}{Y} . \end{equation*}
\end{lemmaIwaniec}
\noindent Since $g ( Xy )$ will be supported on $y \asymp 1 /
X$ and $f ( z )$ is bounded, we have
\begin{eqnarray}
  \nonumber I_f ( X ) & \ll_{f, g} & \int_{y \asymp 1 / X} \int^{1 / 2}_{- 1 / 2} |u ( z
  ) |^2 d \mu z\\
  & \ll_{f, g} & X. \label{3.4}
\end{eqnarray}

For an alternate result to (\ref{3.3}), we evaluate $I_f ( X )$ by shifting the line
of integration.  Starting with equation (\ref{3.2}) and moving the contour of
integration to the line $\textnormal{Re}(s) = 1/2$, we pick up the main term
\begin{equation} M_f ( X ) = \frac{3}{\pi} G ( - 1 ) X < fu, u > \label{Mf} \end{equation}
from the pole of the Eisenstein series at $s = 1$ and obtain the bound
$R_f ( X ) \ll X$ from (\ref{3.4}) where $R_f ( X )$ is the remaining term
\begin{equation}
  R_f ( X ) =  \int_{\Gamma \backslash \mathbbm{H}} p(z) f ( z ) |u ( z ) |^2 d
  \mu z 
\end{equation}
with
\begin{equation*}
p(z):=\frac{1}{2 \pi i} \int_{( 1 / 2 )} G ( - s ) X^s E ( z, s ) ds.
\end{equation*}
With this notation, we now write
\begin{equation}
  I_f ( X ) = M_f ( X ) + R_f ( X ) \label{3.5555}
\end{equation}
and would like to have $R_f ( X )$ small.

From the Fourier series expansion (\ref{2.4}) for the Eisenstein series $E ( z, s )$
and (\ref{B.1}) in Appendix B we have
\begin{eqnarray*}
  E ( z, s ) & = & y^s + \varphi ( s ) y^{1 - s} + O ( |s|^2 y^{- 3 / 2} )\\
  & \ll & \sqrt{y} + |s|^2 y^{- 3 / 2}
\end{eqnarray*}
on the line $\textnormal{Re} ( s ) = 1 / 2$.  Hence $p ( z ) \ll \sqrt{yX}$ if $y
\geqslant 1 / 2$.  Assuming that $\sqrt{y} |f ( z ) |$ is bounded on
$\mathbbm{H}$, we conclude that $R_f ( X ) \ll_{f,g} \sqrt{X} .$

Note that this condition is satisfied for $f ( z )$ an
incomplete Eisenstein series or a cusp form, and so
for these cases $R_f ( X ) \ll_{f,g} \sqrt{X}$.  However, even in the simple case of the constant function, we can not make this conclusion. 
Having the two formulas (\ref{3.3}) and
(\ref{3.5555}) for $I_f ( X )$ brings us to
\begin{equation}
 M_f(X) = \int^{\infty}_0 g ( Xy ) y^{- 2} \left( \int^{1 / 2}_{- 1 / 2} f ( z
  ) |u ( z ) |^2 dx \right) dy + O ( \sqrt{X} ) \label{3.1000}
\end{equation}
which, upon recalling the definition of $M_f(X)$ in (\ref{Mf}), provides the relation between $< fu, u >$ and the relevant shifted
convolution sums which come out of the inner integral.  The main term on the
right, will appear from the zero-th coefficient of $f ( z )$ which is given as
\begin{equation} a_0 ( y ) = \frac{3}{\pi} < f, 1 > + O ( X^{- 1 / 2} ) \label{3.6} . \end{equation}
This is consistent with Lemma A.1 since $< f, 1 > =  a_0 ( y
)=0$ for $f$ a cusp form and $g ( Xy )$ is supported on $y \asymp 1 / X$. 

Now that we've unfolded the integral, our plan is to start picking away at 
parts of the remaining integral in (\ref{3.1000}) which we will be able to show are small until we are left with a shifted convolution sum in the appropriate range. Writing $f^{\ast} ( z ) : = f ( z ) - a_0 ( y )$ and substituting into
(\ref{3.1000}) we get
\begin{equation*}
  M_f(X) = \{ \frac{3}{\pi} < f, 1 > + O (
  X^{- 1 / 2} ) \} I_1 ( X ) + I_{f^{\ast}} ( X ) + O ( \sqrt{X} )
\end{equation*}
where
\begin{eqnarray}
  I_{f^{\ast}} ( X ) & = & \int^{\infty}_0 g ( Xy ) y^{- 2} \left( \int^{1 /
  2}_{- 1 / 2} f^{\ast} ( z ) |u ( z ) |^2 dx \right) dy, \textnormal{ and} \label{3.13}\\
  \nonumber I_1 ( X ) & = & \frac{1}{2 \pi i} \int_{( \sigma )} G ( - s ) X^s
  \int_{\Gamma \backslash \mathbbm{H}} E ( z, s ) |u ( z ) |^2 d \mu zds,  \textnormal{ for }  \sigma > 1.
\end{eqnarray}
Hence, upon moving $\sigma$ to the left in $I_1$,
\begin{equation*}
  M_f(X) = \{ \frac{3}{\pi} < f, 1 > + O (
  X^{- 1 / 2} ) \} \{ \frac{3}{\pi} G ( - 1 ) X + R_1 ( X ) \} + I_{f^{\ast}} ( X ) + O ( \sqrt{X} )
\end{equation*}
where
\begin{equation*} R_1 ( X ) = \int_{\Gamma \backslash \mathbbm{H}} p ( z ) |u ( z ) |^2 d \mu z. \end{equation*}
We see that the pole at $s = 1$ contributes the correct main term
\begin{equation*} \frac{3}{\pi} < f, 1 > \frac{3}{\pi} G ( - 1 ) X \end{equation*}
and we can write
\begin{equation} M_f(X) = \frac{3}{\pi} < f, 1 > \frac{3}{\pi}
   G ( - 1 ) X + I_{f^{\ast}} ( X ) + O ( \sqrt{X} ) \label{3.12}
\end{equation}
provided $R_1 ( X ) \ll \sqrt{X}$.  However, $\sqrt{y}$ is
unbounded on $\mathbbm{H}$ and therefore we have to apply some other analysis
to properly bound $R_1 ( X )$.  We appeal to the Rankin-Selberg zeta function
\begin{equation*} Z ( s ) = \sum_{n \neq 0} | \rho ( n ) |^2 |n|^{- s} \end{equation*}
at the cost of introducing the spectral parameter coming from analysis of the
$K$-Bessel function.

We write
\begin{eqnarray*}
  R_1 ( X ) & = & \frac{1}{2 \pi i} \int_{( 1 / 2 )} G ( - s ) X^s 
  \int_{\Gamma \backslash \mathbbm{H}} E ( z, s ) |u ( z ) |^2 d \mu z ds\\
  & = & \frac{1}{2 \pi i} \int_{( 1 / 2 )} G ( - s ) X^s  \sum_{n \neq 0} |
  \rho ( n ) |^2  \left( \int_0^{\infty} y^{s - 1} K^2_{ir} ( 2 \pi |n|y ) dy
  \right) ds\\
  & = & \frac{1}{2 \pi i} \int_{( 1 / 2 )} G ( - s ) \left( \frac{X}{2 \pi}
  \right)^s Z ( s ) \left( \int_0^{\infty} y^{s - 1} K^2_{ir} ( y ) dy \right)
  ds\\
 & \ll & \frac{X^{1/2}}{r^{1/2}} e^{-\pi r} \int_{( 1 / 2 )} |G ( - s ) Z ( s )| |ds|
\end{eqnarray*}
by (\ref{A.12}) in Lemma B.2. Therefore, the remainder $R_1 ( X )$ will contribute $O ( \sqrt{X} )$, independent of the spectral parameter $r$,
provided $Z ( s )$ satisfies the bound
\begin{equation} Z ( s ) \ll |s|^{10} r^\frac{1}{2} e^{\pi r}, \textnormal{ for Re}(s)=1/2 . \label{logfreeconvex}
 \end{equation}
Assuming this ``log-free" convexity bound, the remaining work will be to simultaneously
evaluate $I_{f^{\ast}} ( X )$ for the cases of $f$ an incomplete Eisenstein
series or a cusp form.  Note that $f^{\ast} ( z )$ is bounded on $\mathbbm{H}$
because both $f ( z )$ and $a_0 ( y )$ are bounded.  We've arrived at the
following lemma.
\begin{lemmaA.2b}
  Let $X \geqslant 1$.  For any $g ( y ) \in C^{\infty}_c (\mathbbm{R}^+ )$,
  any Hecke-Maass cusp form $u(z)$, and any Hecke-Maass cusp form or incomplete Eisenstein series $f(z)$
  with Fourier expansion
  \begin{equation*} f ( z ) = a_0(y)+f^{\ast}(z) = a_0(y)+\sum_{\ell \neq 0} a_{\ell} ( y ) e ( \ell x ), \end{equation*}
  define the integral
  \begin{equation*} I_{f^{\ast}} ( X ) : = \int^{\infty}_0 g ( Xy ) y^{- 2} \left( \int^{1 /
     2}_{- 1 / 2} f^{\ast} ( z ) |u ( z ) |^2 dx \right) dy. \end{equation*}  
  Then for $f(z)$ a cusp form,
  \begin{equation*}\frac{3}{\pi} G ( - 1 ) X < fu, u > = I_{f^{\ast}} ( X ) + O ( X^{1 /
     2} ) . 
\end{equation*}
  Furthermore, if the Rankin-Selberg zeta function associated with $u(z)$
\begin{equation*}
Z ( s ) := \sum_{n \neq 0} | \rho ( n ) |^2 |n|^{- s}
\end{equation*}
 satisfies the bound
 \begin{equation*} Z(s) \ll |s|^{10} r^\frac{1}{2} e^{\pi r} \textnormal{ for Re}(s)=1/2,  \end{equation*}
then
  \begin{equation*}\frac{3}{\pi} G ( - 1 ) X < fu, u > =\frac{3}{\pi} G ( - 1 ) X \frac{3}{\pi} < f, 1 > + I_{f^{\ast}} ( X ) + O ( X^{1 / 2} ). \end{equation*}
\end{lemmaA.2b}
\noindent Our goal now is to show that $I_{f^{\ast}} ( X )$ is much smaller than $X$.
\subsection{Truncation of the fixed function}
The fixed function $f^{\ast} ( z )$ from the previous section will have the
Fourier expansion
\begin{equation} f^{\ast} ( z ) = \sum_{\ell \neq 0} a_{\ell} ( y ) e ( \ell x ) \label{4.111111111} \end{equation}
with coefficient bounds as seen in Lemma A.1.  If $f^{\ast} ( z )$
came from an incomplete Eisenstein series, then we find that the tail of (\ref{4.111111111})
with $| \ell | \geqslant L$ is bounded by
\begin{equation*} y^{- 3 / 2}  \sum_{| \ell | \geqslant L} \frac{\tau ( \ell )}{\ell^2} \ll
   y^{- 3 / 2}  \frac{\log L}{L} \asymp X^{3 / 2}  \frac{\log L}{L} \end{equation*}
for $y \asymp 1 / X$ which is in the range of the integral (\ref{3.13}).  Hence
these terms contribute to $I_{f^{\ast}} ( X )$ at most
\begin{equation*} X^{3 / 2} L^{- 1} ( \log L ) I_1 ( X ) \ll X^{5 / 2} L^{- 1} \log L. \end{equation*}
A similar argument works when $f^{\ast} ( z )$
comes from a fixed cusp form $f ( z )$ in which case $f^{\ast} ( z ) = f ( z
)$.  In this case, we use the bound [K-Sa]
\begin{equation} \rho ( \ell ) \ll \ell^{7 / 64} \tau(\ell) \label{KSa}\end{equation}
to show that the tail contributes at most $X^{5 / 2} {L^{-57 / 64}} \log L$.
Choosing $L = X^2$ in either case gives us
\begin{equation*} I^{\ast}_f ( X ) = I_{f_L} ( X ) + O(X^{3/4}) \end{equation*}
where $I_{f_L} ( X )$ is given by (\ref{3.13}) with $f ( z )$ replaced by
\begin{equation*} f_L ( z ) = \sum_{0 < | \ell | < L} a_{\ell} ( y ) e ( \ell x ) . \end{equation*}
Notice that $f_L ( z )$ is bounded on $\mathbbm{H}$. Now we are left with
showing that
\begin{equation} I_{f_L} ( X ) = \int^{\infty}_0 g ( Xy ) y^{- 2} \left( \int^{1 / 2}_{- 1
   / 2} f_L ( z ) |u ( z ) |^2 dx \right) dy \label{4.3} \end{equation}
has order of magnitude smaller than $X$.

\subsection{Truncation of the variable cusp form}
Recall that a cusp form $u ( z )$ has the Fourier expansion
\begin{equation*} u ( z ) = \sqrt{y}  \sum_{n \neq 0} \rho ( n ) K_{ir} ( 2 \pi |n|y ) e ( nx
   ) . \end{equation*}
We shall see that the key contribution to $I_{f_L} ( X )$ in (\ref{4.3}) comes from
the middle coefficients of $u ( z )$ with $|n|$ about $rX$ for $y \asymp 1 /
X$.  Therefore we write
\begin{equation*} u ( z ) = u^{\ast} ( z ) + u_1 ( z ) + u_2 ( z ) \end{equation*}
where
\begin{equation*} u^{\ast} ( z ) = \sqrt{y}  \sum_{n \neq 0} h \left( \bigg| \frac{n}{rX} \bigg|
   \right) \rho ( n ) K_{ir} ( 2 \pi |n|y ) e ( nx ) \end{equation*}
and $h ( v )$ is some smooth, compactly supported function with $h(v)=1$ for $v\in [1/V,V]$ and $h(v)=0$ for $v$ outside of $[1/2V,2V]$, with $1 \leqslant V \leqslant X$ to be chosen appropriately.  Then $u_1 ( z )$
is the missing partial sum with $|n| \ll rX / V$ and $u_2 ( z )$ is the
missing partial sum with $|n| \gg rXV$.  Put $u_3 ( z ) = u_1 ( z ) + u_2 ( z
)$, so that
\begin{eqnarray*}
  |u^{\ast} ( z ) |^2 & = & |u ( z ) - u_3 ( z ) |^2\\
  & = & |u ( z ) |^2 - 2 \textnormal{Re}\left(u ( z ) u_3 ( z )\right) + |u_3 ( z ) |^2 .
\end{eqnarray*}
Hence, by the boundedness of $f_L ( z )$ we derive from (\ref{4.3}) that
\begin{eqnarray*}
  I_{f_L} ( X ) & = & \int^{\infty}_0 g ( Xy ) y^{- 2} \left( \int^{1 / 2}_{-
  1 / 2} f_L ( z ) |u^{\ast} ( z ) |^2 dx \right) dy\\
  &  & + O \left( \int^{\infty}_0 \left( \int^{1 / 2}_{- 1 / 2} ( |u ( z )
  u_3 ( z ) | + |u_3 ( z ) |^2 ) dx \right) g ( Xy ) y^{- 2} dy \right)\\
  & = & I^{\ast}_{f_L} ( X ) + O ( R ( X ) ),
\end{eqnarray*}
say.  By the Cauchy-Schwartz inequality we get
\begin{equation} R ( X ) \leqslant I_1 ( X )^{1 / 2} R_3 ( X )^{1 / 2} + R_3 ( X ) \label{5.5}
\end{equation}
where
\begin{eqnarray*}
  I_1 ( X ) & = & \int^{\infty}_0 \left( \int^{1 / 2}_{- 1 / 2} |u ( z )
  |^2 dx \right) g ( Xy ) y^{- 2} dy \\
  R_3 ( X ) & = & \int^{\infty}_0 \left( \int^{1 / 2}_{- 1 / 2} |u_3 ( z )
  |^2 dx \right) g ( Xy ) y^{- 2} dy.
\end{eqnarray*}
We already know that $I_1 ( X ) \ll X$, so it remains to estimate $R_3 ( X )$.
We have
\begin{eqnarray*}
  R_3 ( X ) & \leqslant & 2 \int^{\infty}_0 \left( \int^{1 / 2}_{- 1 / 2}
  |u_1 ( z ) |^2 dx \right) g ( Xy ) y^{- 2} dy\\
 & & + 2 \int^{\infty}_0 \left(
  \int^{1 / 2}_{- 1 / 2} |u_2 ( z ) |^2 dx \right) g ( Xy ) y^{- 2} dy\\
  & = & 2 R_1 ( X ) + 2 R_2 ( X )
\end{eqnarray*}
and we will estimate $R_1 ( X )$ and $R_2 ( X )$ directly.  We have
\begin{eqnarray*}
  R_1 ( X ) & \leqslant & \sum_{|n| < rX / V} | \rho ( n ) |^2 \int_0^{\infty}
  g ( Xy ) K^2_{ir} ( 2 \pi |n|y ) y^{- 1} dy\\
  R_2 ( X ) & \leqslant & \sum_{|n| > rXV} | \rho ( n ) |^2 \int_0^{\infty} g
  ( Xy ) K^2_{ir} ( 2 \pi |n|y ) y^{- 1} dy
\end{eqnarray*}
and by Theorem B.5 and Corollary B.2,
\begin{eqnarray*}
  R_1 ( X ) & \ll & \sum_{|n| < rX / V} | \rho ( n ) |^2 \frac{1}{r} e^{-\pi r}\\
  R_2 ( X ) & \ll & \sum_{|n| > rXV} | \rho ( n ) |^2 \left(\frac{X}{n}\right)^{2} r e^{-\pi r}
\end{eqnarray*}
so that
\begin{eqnarray*}
  R_1 ( X ) \ll \frac{X}{V} & \textnormal{and} & R_2 ( X ) \ll \frac{X}{V} .
\end{eqnarray*}
Hence $R_3 ( X ) \ll X / V$ and by (\ref{5.5}) we conclude that
\begin{equation*} R ( X ) \ll XV^{- 1 / 2} . \end{equation*}
This bound is smaller than the main term in (\ref{3.12}) if $V$ is chosen to be
sufficiently large. 
\subsection{Opening $I^{\ast}_{f_L} ( X )$}
We are left with evaluating
\begin{equation} I^{\ast}_{f_L} ( X ) = \int^{\infty}_0 g ( Xy ) y^{- 2} \left( \int^{1 /
   2}_{- 1 / 2} f_L ( z ) |u^{\ast} ( z ) |^2 dx \right) dy \label{6.1} \end{equation}
where
\begin{equation*} f_L ( z ) = \sum_{0 < | \ell | < L} a_{\ell} ( y ) e ( \ell x ) \end{equation*}
and
\begin{equation*} u^{\ast} ( z ) = \sqrt{y}  \sum_{n \neq 0} h \left( \bigg| \frac{n}{rX} \bigg|
   \right) \rho ( n ) K_{ir} ( 2 \pi |n|y ) e ( nx ) . \end{equation*}
Recall that $h ( v )$ is a nice compactly
supported function at our disposal and that $L$ can be taken to be $X^2$. 
Hence
\begin{equation} I^{\ast}_{f_L} ( X ) = \sum_{0 < |\ell| \leqslant L} S_{\ell} ( X ) \label{4.2}
\end{equation}
where
\begin{eqnarray}
 \nonumber S_{\ell} ( X ) & = & \sum_n h \left( \bigg| \frac{n}{rX} \bigg| \right) h \left( \bigg|
  \frac{n + \ell}{rX} \bigg| \right) \overline{\rho ( n )} \rho ( n + \ell )\\
  &  & \times \int^{\infty}_0 a_{\ell} ( y ) g ( Xy ) K_{ir} ( 2 \pi |n|y )
  K_{ir} ( 2 \pi |n + \ell |y ) y^{- 1} dy. \label{6.3}
\end{eqnarray}
From Lemma A.2b, \S A.3, \S A.4 and (\ref{6.3}) we can now state the following theorem.
\begin{theoremA.5}
  Let $X \geqslant 1$, $S_{\ell} ( X )$ be as in \textnormal{(\ref{6.3})}, and let $u(z)$ be the variable cusp form with Fourier expansion
\begin{equation*}
u ( z ) = \sqrt{y}  \sum_{n \neq 0} \rho ( n ) K_{ir} ( 2 \pi |n|y ) e ( nx ).
\end{equation*}
Then for any fixed cusp form $f ( z )$, 
  \begin{equation*} < fu, u > = \frac{\pi}{3} \left(G( - 1 ) X\right)^{-1} \sum_{0 < | \ell | \leqslant
     X^2} S_{\ell} ( X ) + O(X^{-\delta}) 
\end{equation*}
  for some constant $\delta>0$. Furthermore, if the
  Rankin-Selberg zeta function associated with $u(z)$
\begin{equation*}
Z ( s ) := \sum_{n \neq 0} | \rho ( n ) |^2 |n|^{- s}
\end{equation*}
 satisfies the bound
  \begin{equation*} Z ( s ) \ll |s|^{10} r^\frac{1}{2} e^{\pi r}, \textnormal{ for Re}(s)=1/2, \end{equation*}
  then for any fixed incomplete Eisenstein series $f(z)$,
  \begin{equation*} < fu, u > = \frac{3}{\pi} < f, 1 > + 
\frac{\pi}{3} \left(G( - 1 ) X\right)^{-1} \sum_{0 < | \ell | \leqslant X^2} S_{\ell} ( X ) + O(X^{-\delta})  
\end{equation*}
for the same constant $\delta>0$.
\end{theoremA.5}

By Theorem A.5, in order to prove (\ref{1}) and QUE, it only remains to show that for each fixed $0<|\ell|<X^2$ and test function $g ( y )$ we have
\begin{eqnarray}
  S_{\ell} ( X ) = o ( 1 ) & \textnormal{as} & r \longrightarrow \infty. \label{rateofconv}
\end{eqnarray}
One may choose to proceed by looking for cancellations in $S_\ell(X)$, but we believe that summing over absolute values should be sufficient for obtaining (\ref{rateofconv}). Appealing to Theorem B.5, we get a convenient normalization of our coefficients $\rho(n)$.  In particular, with each $a_\ell(y)$ bounded as in Lemma A.1, we can bound $S_\ell(X)$ by 
\begin{eqnarray}
S_\ell(X)& \ll & \sum_{n\asymp r} |\overline{\rho(n)} \rho(n+\ell)| \int_{0}^{\infty} g(y) K_{ir}^2\left(\frac{2\pi|n|}{X} y \right)y^{-1} dy \nonumber \\ 
& \ll &\sum_{n\asymp r}  \frac{|\overline{\rho ( n )} \rho ( n + \ell )|}{r \cosh( \pi r )} = \sum_{n\asymp r}  \frac{|\lambda ( n ) \lambda ( n + \ell )|}{r L(\textnormal{sym}^2 u,1)} \nonumber
\end{eqnarray}
with the implied constant depending on $\ell$, the test functions $h,g$ and $X$. Therefore, if one shows that for fixed $0<|\ell|<X^2$,
\begin{equation*}
  \sum_{n\asymp r} \left|{\lambda} ( n ) \lambda ( n + \ell )\right| = o ( r L(\textnormal{sym}^2 u, 1) )
\end{equation*}
as $r\longrightarrow \infty$, then applying Theorem A.5 and sending $X\longrightarrow \infty$ arbitrarily slowly will prove (\ref{1}) and QUE modulo the ``log-free" convexity bound (\ref{logfreeconvex}).

\newpage
\section{Properties of the $K$-Bessel function}
The modified Bessel functions are solutions to the second order differential
equation
\begin{equation*} z^2 f'' + zf' - ( z^2 + \nu^2 ) f = 0. \end{equation*}
One solution for $f$ is given by the power series
\begin{equation*} I_{\nu} ( z ) : = \sum^{\infty}_{k = 0} \frac{1}{k! \Gamma ( k + 1 + \nu )}
   \left( \frac{z}{2} \right)^{\nu + 2 k} . \end{equation*}
A second solution is given by
\begin{equation} K_{\nu} ( z ) : = \frac{\pi}{2} ( \sin \pi \nu )^{- 1} ( I_{- \nu} ( z ) -
   I_{\nu} ( z ) ). \label{A.1} \end{equation}
We refer to $\nu$ as
the ``order'' and $z$ as the ``argument'' of the Bessel function.

In the real variable $y$, the Bessel function $K_{\nu} ( y )$ can be well
approximated outside of the range $y \asymp 1 + | \nu |$.  We can get a good
approximation for small $y$ from the power series expansion and have ([Iw],
B.36)
\begin{equation}
  K_{\nu} ( y ) = \left( \frac{\pi}{2 y} \right)^{1 / 2} e^{- y}  \left( 1 + O
  \left( \frac{1 + | \nu |^2}{y} \right) \right) \label{A.2}, 
\end{equation}
which gives exponential decay in $y$ for $y > 1 + | \nu |^2$.

Two useful formulas are the integral representation ([Iw], p. 227)
\begin{equation} K_{\nu} ( z ) = \pi^{- 1 / 2} \Gamma \left( \nu + \frac{1}{2} \right)
   \left( \frac{z}{2} \right)^{- \nu}  \int_0^{+ \infty} ( t^2 + 1 )^{- \nu -
   1 / 2} \cos ( tz ) dt \label{A.3} \end{equation}
and the Mellin transform ([Iw], p. 228)
\begin{equation} \int_0^{+ \infty} K_{\mu} ( y ) K_{\nu} ( y ) y^{s - 1} dy = 2^{s - 3}
   \Gamma ( s )^{- 1} \prod \Gamma \left( \frac{s \pm \mu \pm \nu}{2} \right)
   \label{A.4} \end{equation}
for $\textnormal{Re} ( s ) > | \textnormal{Re} \mu | + | \textnormal{Re} ( \nu ) |$.
From the integral representation (\ref{A.3}), the following bound is known for real $z$ and imaginary order $ir$
  \begin{equation} K_{ir} ( y ) \ll | \Gamma ( \frac{1}{2} + ir ) | \frac{1 + r^2}{y^2} . \label{B.1} \end{equation}
Furthermore, when $\mu = \nu = ir$, with $r \in \mathbbm{R}^+$ in the Mellin transform (\ref{A.4}),
we get the following bound for $s = \sigma + it$ with $\sigma > 0$ by Stirling's formula,
  \begin{equation} \int^{+ \infty}_0 K^2_{\nu} ( y ) y^{s - 1} dy \ll_{\sigma} \frac{1}{(
     |t| + \sigma )^{1 / 2}} \min \left\{ r^{\sigma - 1} e^{- \pi r}, |
     t|^{\sigma - 1} e^{- \pi |t| / 2} \right\} \label{A.11} . \end{equation}
In particular, for $0 < \sigma \leqslant 3 / 2$ we have
\begin{equation} 
\int^{+ \infty}_0 K^2_{\nu} ( y ) y^{s - 1} dy \ll_{\sigma}  r^{\sigma - 1} e^{- \pi r}.  \label{A.12}
\end{equation}

For $s = 0$ we see that the integral can not be bounded in this manner.
However, a bound can be established when integrating against a smooth
compactly supported function $g$ with Mellin transform
\begin{equation*} G ( s ) := \int_0^{+ \infty} g ( y ) y^{s - 1} dy \end{equation*}
which will be entire and satisfies
\begin{equation} G ( s ) \ll_{g, A} ( |s| + 1 )^{- A} \label{A.13} \end{equation}
for any $A \geqslant 0$, uniformly in vertical strips.  Therefore, for $g \in
C^{\infty}_c (\mathbbm{R}^+ )$,
\begin{equation*} g ( y ) = \frac{1}{2 \pi i} \int_{( \sigma )} G ( s ) y^{- s} ds \end{equation*}
converges for all $- B \leqslant \sigma \leqslant B$ for any $B$ with $y > 0$.
The integral now under consideration is
\begin{equation} \int^{+ \infty}_0 g ( y ) K^2_{\nu} ( y ) y^{- 1} dy = \frac{1}{2 \pi i}
   \int_{( \sigma )} G ( - s ) \int^{+ \infty}_0 K^2_{\nu} ( y ) y^{s - 1}
   dyds. \label{A.14} \end{equation}
The contribution of $|t|^{\sigma - 1}$ for large $\sigma$ and large $t$ in
(\ref{A.11}) is now dominated by our function $G ( s )$, provided we take $A$ large
enough in (\ref{A.13}).  Upon integrating with respect to $s$ we arrive at the
following lemma.
\begin{lemmaB.1}
  For any smooth, compactly supported function $g ( y )$, order $\nu = ir$
  with $r > 0$, and $\sigma > 0$ we have
  \begin{equation} \int^{+ \infty}_0 g ( y ) K^2_{\nu} ( y ) y^{- 1} dy \ll_{\sigma, g}
     r^{\sigma - 1} e^{- \pi r} . \label{A.15} \end{equation}
\end{lemmaB.1}
\begin{corollaryB.2}
  For any smooth, compactly supported function $g ( y )$, constant $w > 0$,
  order $\nu = ir$ with $r > 0$, and $\sigma > 0$ we have
  \begin{equation} \int^{+ \infty}_0 g ( y ) K^2_{\nu} ( wy ) y^{- 1} dy \ll_{\sigma, g} 
     \frac{1}{r} e^{- \pi r}  \left( \frac{r}{w} \right)^{\sigma} . \label{A.16}
  \end{equation}
\end{corollaryB.2}

Corollary B.2 is particularly nice for constants $w \geqslant r$. In order to
establish an equally good result for small $w$, we return to the Mellin transform (\ref{A.4}), apply the duplication formula for the gamma function and
move the contour line $( \sigma )$ in (\ref{A.14}) to the left of $s = 0$. However,
$\Gamma ( s )$ has a simple pole at each negative integer so we take care to
not move the contour so far as to pick up any of these poles.  Moving the
contour to any $\textnormal{Re} ( s ) = \delta > - 2$, we pick up the simple
poles of the $\Gamma$ factors at $s = 0, 2 ir, - 2 ir$, and our integral (\ref{A.14}) equals
\begin{equation*}
  \frac{1}{2 \pi i} \int_{( \delta )} G ( - s ) \left( \int^{+\infty}_0 K^2_{\nu} ( y ) y^{s - 1} dy \right) ds+\sum_{s=-2ir,0,2ir} \textnormal{res}_{s} \textnormal{GF} (s)
\end{equation*}
where
\begin{equation*} \textnormal{GF} ( s ) = \frac{\sqrt{\pi}}{4}  \frac{\Gamma \left( \frac{s}{2}
   \right)}{\Gamma \left( \frac{s + 1}{2} \right)} \Gamma \left( \frac{s + 2
   ir}{2} \right) \Gamma \left( \frac{s - 2 ir}{2} \right) G ( - s ) . \end{equation*}
Since the residue at $s = 0$ of $\Gamma ( s )$ is 1 and $\Gamma ( 1 / 2 ) =
\sqrt{\pi}$, the residues of $\textnormal{GF} ( s )$ are
\begin{eqnarray*}
  \textnormal{res}_{s = 0} \textnormal{GF} ( s ) & = & \frac{1}{4} \Gamma ( ir ) \Gamma
  ( - ir ) G ( 0 ),\\
  \textnormal{res}_{s = 2 ir} \textnormal{GF} ( s ) & = & \frac{1}{8} \Gamma^2 ( ir ) G (
  - 2 ir ),\\
  \textnormal{res}_{s = - 2 ir} \textnormal{GF} ( s ) & = & \frac{1}{8} \Gamma^2 ( - ir )
  G ( 2 ir )
\end{eqnarray*}
and because $G ( \pm 2 ir ) \ll ( r + 1 )^{- A}$, the residues contribute
\begin{equation} \frac{1}{4} \Gamma ( ir ) \Gamma ( - ir ) \{ G ( 0 ) + O ( r^{- 2} ) \} .
   \label{A.17} \end{equation}
Integrating on the line $\textnormal{Re} ( s ) = \delta = - 1$ gives
\begin{equation} \frac{1}{2 \pi i} \int_{( - 1 )} G ( - s ) \left( \int^{+ \infty}_0
   K^2_{\nu} ( y ) y^{s - 1} dy \right) ds \ll_g  \frac{1}{r^2} e^{- \pi r}, \label{A.18} \end{equation}
by similar analysis which led us to Lemma B.1, and also gives the following Lemma.
\begin{lemmaB.3}
  For any smooth, compactly supported function $g ( y )$, order $\nu = ir$
  with $r > 0$, we have
  \begin{equation} \int^{+ \infty}_0 g ( y ) K^2_{\nu} ( y ) y^{- 1} dy = \frac{1}{4}
     \Gamma ( ir ) \Gamma ( - ir ) \{ G ( 0 ) + O ( r^{- 2} ) \} + O_g (
     \frac{1}{r^2} e^{- \pi r} ) . \label{A.19} \end{equation}
\end{lemmaB.3}
\noindent Introducing the constant $w$ gives
\begin{corollaryB.4}
  For any smooth, compactly supported function $g ( y )$, constant $w > 0$,
  order $\nu = ir$ with $r > 0$, we have
  \begin{equation*} \int^{+ \infty}_0 g ( y ) K^2_{\nu} ( wy ) y^{- 1} dy = \frac{1}{4}
     \Gamma ( ir ) \Gamma ( - ir ) \{ G ( 0 ) + O ( r^{- 2} ) \} + O_g (
     \frac{w}{r^2} e^{- \pi r} ) .  \end{equation*}
\end{corollaryB.4}
\noindent By Corollary B.4, since $\Gamma ( ir ) \Gamma ( - ir ) \ll \frac{1}{r} e^{- \pi
r}$, we get
\begin{equation} \int^{+ \infty}_0 g ( y ) K^2_{\nu} ( wy ) y^{- 1} dy \ll_g \frac{1}{r}
   e^{- \pi r} + \frac{w}{r^2} e^{- \pi r} . \label{A.20} \end{equation}
Combining this with Corollary B.2, we see that the same bound exists for any attached positive constant
$w$.
\begin{theoremB.5}
  For any smooth, compactly supported function $g ( y )$, constant $w > 0$,
  order $\nu = ir$ with $r > 0$, we have
  \begin{equation} \int^{+ \infty}_0 g ( y ) K^2_{\nu} ( wy ) y^{- 1} dy \ll_g \frac{1}{r}
     e^{- \pi r} . \label{A.21} \end{equation}
\end{theoremB.5}
\noindent This establishes all of the necessary properties for our applications.

\end{document}